\documentclass{article}

\usepackage {amsmath}
\usepackage{color}
\usepackage {hyperref}
\usepackage {amsfonts}
\usepackage {stmaryrd}
\usepackage {amsthm}
\usepackage[flushmargin]{footmisc}
\usepackage[parfill]{parskip}
\usepackage {amssymb}
\usepackage {fullpage}
\usepackage {mathrsfs}
\usepackage{mathtools}
\usepackage {amscd}
\usepackage{bbm}
\usepackage[all,cmtip]{xy}
\DeclareMathAlphabet{\mathpzc}{OT1}{pzc}{m}{it}

\title{Vanishing sheaves and the geometric Whittaker model}
\author{Roman Bezrukavnikov and Tanmay Deshpande}
\date{}

\newtheorem {theorem} {Theorem} [section]
\newtheorem {proposition} [theorem] {Proposition}
\newtheorem {conjecture} [theorem] {Conjecture}
\newtheorem {lemma} [theorem] {Lemma}
\newtheorem {corollary} [theorem] {Corollary}

\theoremstyle{definition}
\newtheorem {defn} [theorem] {Definition}

\newtheorem {prob} [theorem]  {Problem}

\newtheorem {rk} [theorem]  {Remark}
\newtheorem {ex} [theorem] {Example}

\newcommand{\beq}{\begin{equation}}
\newcommand{\eeq}{\end{equation}}
\newcommand{\bthm}{\begin {theorem}}
\newcommand{\ethm}{\end {theorem}}
\newcommand{\bprop}{\begin {proposition}}
\newcommand{\eprop}{\end {proposition}}
\newcommand{\bprob}{\begin {prob}}
\newcommand{\eprob}{\end {prob}}
\newcommand{\bcor}{\begin {corollary}}
\newcommand{\ecor}{\end {corollary}}
\newcommand{\blem}{\begin{lemma}}
\newcommand{\elem}{\end{lemma}}
\newcommand{\bdefn}{\begin{defn}}
\newcommand{\edefn}{\end{defn}}
\newcommand{\bconj}{\begin{conjecture}}
\newcommand{\econj}{\end{conjecture}}
\newcommand{\brk}{\begin{rk}}
\newcommand{\erk}{\end{rk}}
\newcommand{\bpf}{\begin{proof}}
\newcommand{\epf}{\end{proof}}
\newcommand{\bex}{\begin{ex}}
\newcommand{\eex}{\end{ex}}
\newcommand{\bit}{\begin{itemize}}
\newcommand{\eit}{\end{itemize}}

\newcommand{\echi}{\mathbbm{e}_\chi}
\newcommand{\ecircchi}{\mathbbm{e}^\circ_\chi}
\newcommand{\bbecirc}{\mathbbm{e}^\circ}
\newcommand{\bbe}{\mathbbm{e}}

\renewcommand {\bar} {\overline}
\renewcommand{\subset}{\subseteq}
\renewcommand{\supset}{\supseteq}

\newcommand{\xto}{\xrightarrow}
\newcommand{\rar}[1]{\stackrel{#1}{\longrightarrow}}
\newcommand{\onto}{\twoheadrightarrow}

\newcommand{\opname}{\operatorname}
\newcommand{\Perv}{\opname{Perv}}
\newcommand{\IC}{\opname{IC}}
\newcommand{\av}{{\opname{av}}}
\newcommand{\Av}{{\opname{Av}}}

\newcommand{\forg}{\opname{Fg}}
\newcommand{\op}{{\opname{op}}}

\newcommand{\ind}{\operatorname{ind}}

\newcommand{\res}{\operatorname{res}}

\newcommand{\ad}{\operatorname{ad}}
\newcommand{\ab}{\operatorname{ab}}
\newcommand{\Hom} {\opname{Hom}}
\newcommand{\Coh} {\opname{Coh}}
\newcommand{\RHom} {\opname{RHom}}

\renewcommand{\Vec}{\opname{Vec}}

\newcommand{\pt}{\opname{pt}}
\newcommand{\id}{\operatorname{id}}
\newcommand{\Id}{\operatorname{Id}}

\newcommand{\Rep}{\opname{Rep}}

\newcommand{\tr}{\opname{tr}}

\newcommand{\reg}{{\opname{reg}}}
\newcommand{\HC}{\opname{HC}}

\newcommand{\bop}{{\overline{\mathfrak{B}}}}
\newcommand{\Dcent}{\D^{\circ}}
\newcommand{\Dcirc}{\D^\circ}
\newcommand{\Wh}{\eL\D(G)\eL}
\newcommand{\Yok}{\eUop\D(G)\eUop}

\newcommand{\bbA}{\mathbb{A}}

\newcommand{\bbZ}{\mathbb{Z}}
\newcommand{\Fp} {\mathbb{F}_p}

\newcommand{\Qlcl} {\overline{\mathbb{Q}}_{\ell}}

\newcommand{\Ga} {\mathbb{G}_{\opname{a}}}
\newcommand{\Gm} {\mathbb{G}_{\opname{m}}}

\newcommand{\Qcal} {\mathcal{Q}}

\renewcommand{\>}{\rangle}

\renewcommand{\u} {\underline}

\renewcommand{\L}{\mathcal{L}}
\newcommand{\chiTzeta}{\chi_{\bar{T}\zeta}}
\newcommand{\D} {\mathscr{D}}

\newcommand{\DGG}{\D_G(G)}

\newcommand{\Gcal}{\mathcal{G}}

\newcommand{\Ccal}{\mathcal{C}}
\newcommand{\Fcal}{\mathcal{F}}

\newcommand{\Kcal}{\mathcal{K}}
\newcommand{\Mcal}{\mathcal{M}}
\newcommand{\Ucal}{\mathcal{U}}
\newcommand{\Zcal}{\mathcal{Z}}
\newcommand{\Ocal}{\mathcal{O}}
\newcommand{\Ecal}{\mathcal{E}}

\newcommand{\Iscr}{\mathscr{I}}
\newcommand{\Yscr}{\mathscr{Y}}

\newcommand{\Pcal}{\mathcal{P}}

\newcommand{\eL}{e_\L}
\newcommand{\eU}{e_U}
\newcommand{\eUop}{e_{\bar{U}}}

\newcommand{\Uop}{\bar{U}}
\newcommand{\Bop}{\bar{B}}
\newcommand{\Pop}{\bar{P}}

\newcommand{\Spr}{\opname{Spr}}
\newcommand{\Irr}{\opname{Irr}}

\newcommand{\Hecke}{\eU\D_B(G)\eU}

\begin{document}
\maketitle
\begin{abstract}
Let $G$ be a connected reductive algebraic group over an algebraically closed field  $k$ of characteristic $p>0$ and let $\ell$ be a prime number different from $p$. Let $U\subset G$ be a maximal unipotent subgroup, and let $T$ be a maximal torus normalizing $U$ with normalizer $N=N_G(T)$. Let $W=N/T$ be the Weyl group of $G$. Let $\L$ be a non-degenerate multiplicative $\Qlcl$-local system on $U$. In this paper we prove that the bi-Whittaker category, namely the triangulated monoidal category of $(U,\L)$-bi-equivariant $\Qlcl$-complexes on $G$, is monoidally equivalent to an explicit thick triangulated monoidal subcategory $\Dcirc_W(T)\subset \D_W(T)$ of ``$W$-equivariant central sheaves'' on the torus, answering a question raised by Drinfeld. In particular, the bi-Whittaker category has the structure of a symmetric monoidal category. We also study a certain thick triangulated monoidal subcategory $\Dcirc_G(G)\subset \D_G(G)$ of ``vanishing sheaves'' and prove that it is braided monoidally equivalent to an explicit thick triangulated monoidal subcategory $\Dcirc_N(T)\subset \D_N(T)$ of ``$N$-equivariant central sheaves'' on the torus. The above equivalence is given by an enhancement of the parabolic restriction functor restricted to the subcategory $\Dcirc_G(G)$. 
\end{abstract}

\setcounter{tocdepth}{1}
\tableofcontents
\section{Introduction and main results}
Let $G$ be a connected reductive algebraic group over an algebraically closed field $k$ of characteristic $p$, which will be assumed to be positive unless stated otherwise. Let us fix a prime number $\ell$ which is invertible in $k$. For a $k$-scheme $X$, let $\D(X)$ denote the $\Qlcl$-linear triangulated category of bounded constructible $\Qlcl$-complexes on $X$. We let $\omega_X\in \D(X)$ denote the dualizing complex on $X$. We let $\D_G(G)$ denote the conjugation equivariant bounded derived category of $\Qlcl$-complexes on $G$. It is a $\Qlcl$-linear triangulated braided monoidal category, with the monoidal structure being defined by convolution with compact support. We refer to \cite{BoDr:14} for more on these structures on $\D_G(G)$.

We will use the following notation related to the connected reductive group $G$. Let $T\subset B$ be a maximal torus and a Borel subgroup of $G$ containing it. Let $U$ be the unipotent radical of $B$ so that $B=TU$. Let $\Bop=T\Uop$ be the opposite Borel. Let $Z\subset T$ be the center of $G$. Let $\Phi\supset\Phi^+\supset\Delta$ denote the set of roots, positive roots and simple roots respectively, which are determined by the choice of the pair $T\subset B$. Let $N=N(T)$ denote the normalizer of the maximal torus $T$ and let $W=N(T)/T$ denote the Weyl group of $G$. We denote the lattice of characters of $T$ by $X^*(T)$ and the dual lattice of cocharacters of $T$ by $X_*(T)$.

Let $\L$ be a non-degenerate multiplicative local system on the maximal unipotent $U$. Note that by a non-degenerate multiplicative local system $\L$ on $U$, we mean a multiplicative local system $\L$ such that $\L|_{U_\alpha}\ncong {\Qlcl}_{U_\alpha}$ for every simple root $\alpha\in \Delta$ and such that $\L|_{U_\beta}\cong {\Qlcl}_{U_\beta}$ for every non-simple positive root $\beta\in \Phi^+\setminus\Delta$. Such non-degenerate multiplicative local systems on $U$ form a single orbit under the adjoint action of $T$ on $U$. Moreover, the action of the adjoint torus $T/Z$ is simply transitive on the space of non-degenerate multiplicative local systems on $U$.

The main goal of this paper is to study the {\em bi-Whittaker} category of $G$, namely the category of $(U,\L)$-bi-equivariant sheaves on $G$, and to identify it with a certain full thick triangulated monoidal subcategory $\Dcirc_W(T)$ of the triangulated symmetric monoidal category $\D_W(T)$ of $W$-equivariant complexes on the maximal torus $T$,  known as the subcategory of {\em central sheaves} on the torus. This answers a question raised by Drinfeld. Analogues of these questions in the setting of $D$-modules have been studied by Ginzburg, Lonergan, Ben-Zvi - Gunningham and Gannon in \cite{Gi:18}, \cite{Lo:18}, \cite{BZGu} and \cite{Gan} respectively. In particular, \cite[Thm. 1.4]{Gan} building on \cite{Gi:18, Lo:18} proves an analogue
of our second main result (Theorem \ref{thm:main2}) in the $D$-module setting. (The results of \cite{Gan}
and of the present paper were obtained independently at about the same time).
In the present paper we work in the context of $\ell$-adic sheaves on schemes over a field of characteristic $p>0$,
however, the arguments also apply to the setting of $D$-modules yielding results close to some of the results in the above mentioned papers. In contrast, their methods are specific to the $D$-modules setting.

As a consequence of the above result, we see that the bi-Whittaker category has the structure of a {\em symmetric} monoidal category. 
This can be thought of as a geometric analogue of the ``uniqueness of the Whittaker model'' or the multiplicity-freeness of the Gelfand-Graev representations.

The category of central sheaves on the torus is also studied by Chen in \cite{Chen:1902,Chen:1909}. The categories that appear above are also closely related to a certain full subcategory $\Dcirc_G(G)\subset \D_G(G)$ of {\em ``vanishing complexes''} on $G$. Before turning to the bi-Whittaker category, we will study this category of vanishing complexes. By definition, the Harish-Chandra restriction agrees with the parabolic restriction functor on this subcategory. In our first main result (Theorem \ref{thm:main1}), we prove that an upgraded version of the parabolic restriction functor defines a triangulated braided monoidal equivalence $\res^\circ:\Dcirc_G(G)\xto{\cong}\Dcirc_N(T)$, where the latter category is the thick triangulated monoidal  subcategory of ``central sheaves'' in $\D_N(T)$. %In the course of proving this result, we will also prove the ``vanishing conjecture'' from \cite{Chen:1902}. 

Whittaker model plays a central role in several chapters of geometric representation theory. In particular, it provides a conceptual framework for Soergel's description of category $O$ \cite{Soe} and its generalization to quasi-split real groups \cite{BeVi} elucidated by our present result
(see Remark \ref{Soe_rem} below). Another potential application is to the study of Iwahori-Whittaker sheaves on the affine flag variety.

Notice also that, while not equivalent, the bi-Whittaker and the vanishing sheaves categories turn out to be closely related; in particular, our results show that the hearts of the perverse $t$-structure on the two categories are equivalent (see Corollary \ref{equiv_ab}).
The vanishing sheaves category can be viewed as the categorification of
an analogue of the (stable) Bernstein center in the context of finite Chevalley groups (see Remarks \ref{cent_an}, \ref{gamma_sh}). Thus the relation between the bi-Whittaker category 
and $\Dcirc_G(G)$ can be viewed as a finite field analogue of the role
of Whittaker model in the local geometric Langlands duality.

In our proof of these results, the Yokonuma-Hecke category $\D(\Uop\backslash G /\Uop)$ will play an important role. In particular, we recall the definition and properties of certain shifted perverse sheaves $\Kcal_w\in \D(\Uop\backslash G/\Uop)$ parametrized by $w\in W$ defined by Kazhdan and Laumon in \cite{KazLau}. These perverse sheaves are geometric analogues of the generators of the Yokonuma-Hecke algebra studied by Juyumaya in \cite{Ju:98}. We describe our main results in more detail below.

\subsection{The bi-Whittaker category and some other unipotent Hecke categories}
Let us begin by introducing our main object of study in this paper, namely the category of bi-Whittaker sheaves on the group $G$. We also introduce some other unipotent Hecke categories. For this it will be convenient to recall some facts and notation from the theory of character sheaves on unipotent groups developed by Boyarchenko-Drinfeld in \cite{BoDr:14}. 
For a multiplicative local system $\L$ on any unipotent algebraic group $H$, let $\eL:=\L\otimes\omega_H=\L[2\dim H](\dim H)$ be the corresponding closed idempotent (see \cite{BoDr:14} for details) in $\D_H(H)$. If $\L$ is the trivial local system on $H$, we will often denote the corresponding idempotent simply by $e_H\in \D_H(H)$. For a connected unipotent group $H$, we have its Serre dual $H^*$ which is a (possibly disconnected) perfect commutative unipotent group over $k$ which parametrizes the multiplicative local systems on $H$.

Let us fix a non-trivial additive character $\psi:\Fp\to \Qlcl^\times$. This defines for us the Artin-Schreier local system $\L_{\psi}$ on $\Ga$ and gives us an identification of the Serre dual $\Ga^*\cong \Ga$ (after passing to the perfectizations, see \cite{BoDr:14} for details), i.e. the moduli space of multiplicative local systems on $\Ga$ gets identified with the perfectization of $\Ga$ once we fix the non-trivial additive character $\psi$. To define the bi-Whittaker category we will fix a non-degenerate multiplicative local system $\L$ on the maximal unipotent subgroup $U\subset G$, as has been described above. 

\brk
The choice of such a non-degenerate multiplicative local system $\L$ on $U$ is equivalent to the choice of a pinning of the reductive group $G$. Namely, given such a local system $\L$ and a simple root $\alpha$ we make the choice of the root homomorphism $f_\alpha:SL_2\to G$ such that $\L$ pulls back to $\L_\psi$ on $\Ga\cong U^+\subset SL_2$ via the isomorphism $f_\alpha|_{U^+}:U^+\xto{\cong} U_\alpha\subset U$, where $\Ga\cong U^+\subset SL_2$ is the subgroup of upper triangular unipotent matrices. The root homomorphism $f_\alpha$ is uniquely determined by this requirement for each simple root $\alpha$, which hence determines the pinning.
\erk

The non-degenerate multiplicative local system $\L$ gives rise to the closed idempotent $\eL\in \D_U(G)$.
\bdefn
(i) The bi-Whittaker category is defined to be the triangulated monoidal category 
$$\Wh\subset \D_U(G)\subset \D(G),$$ i.e. it is the Hecke subcategory corresponding to the closed idempotent $\eL$ in the terminology of \cite{BoDr:14}. The object $\eL\in \Wh$ is the unit object.\\
(ii) The Yokonuma-Hecke category is defined to be the triangulated monoidal category $\D(U\backslash G/U)=e_U\D(G)e_U\subset \D_U(G)\subset \D(G)$ whose unit object is $e_U$. In this paper it will be more convenient for us to consider the Yokonuma-Hecke category $\D(\Uop\backslash G/\Uop)=\eUop\D(G)\eUop$, with respect to the opposite unipotent subgroup, whose unit object is $\eUop$.\\
(iii) Note that we can consider the closed idempotent $e_U$ as a closed idempotent in $\D_B(G)$. By the Hecke category with respect to $U$ we mean the triangulated monoidal category $e_U\D_B(G)e_U=e_U\D_B(G)=\D_B(G)e_U$ with unit object $e_U$. Sometimes it will be more convenient to consider the Hecke category $\eUop\D_{\Bop}(G)\eUop$ with respect to the opposite unipotent $\Uop$.
\edefn

\subsection{The Harish-Chandra functors and their adjoints}
Let us now describe some functors which relate the above Hecke categories to the category $\D_G(G)$. Recall that for a finite group, the convolution algebra of class functions is the center of the convolution algebra of all functions on the finite group. The geometric analogue of this is the fact that the forgetful functor $\D_G(G)\to \D(G)$ is central, i.e. it naturally factors through the Drinfeld center $\Zcal(\D(G))\to \D(G)$ (see \cite{BoDr:14}). In particular we have the functorial braiding isomorphisms $\beta_{A,C}:A\ast C\xto{\cong}C\ast A$ in $\D(G)$ for $A\in \D_G(G), C\in \D(G)$.

Using the above, we can define the central monoidal functors 
\beq
\HC:\DGG\to\eU\D_{B}(G)\eU\subset \D_{B}(G) \mbox{ defined by } A\mapsto \eU\ast A\cong A\ast \eU,
\eeq
\beq
\HC_\L:\DGG\to \eL\D(G)\eL\subset \D_U(G) \mbox{ defined by } A\mapsto \eL\ast A\cong A\ast \eL.
\eeq 

Note that since the flag variety $G/B$ is projective,  the canonical morphism between the averaging functors with and without compact supports is an isomorphism $\av_B^G\rar{\cong}\Av_B^G$. We refer to \cite{BoDr:14, Desh:16} for more about the averaging functors. Since $\Av_{(\cdot)}^G$ is right adjoint to the forgetful functor, the averaging functors define right adjoints to the Harish-Chandra functors defined above:
\beq
\av_{B}^{G}=\Av_{B}^G:\eU\D_{B}(G)\eU\subset \D_{B}(G)\to \D_G(G),
\eeq
\beq
\Av_U^G:\eL\D(G)\eL\subset \D_U(G)\to \D_G(G).
\eeq
Hence we have the two adjoint pairs
\beq\label{eq:hcavadjunction}
(\HC,\Av_B^G=\av_B^G) \mbox{ and } (\HC_\L,\Av_U^G).
\eeq
Note that for $A\in \D_G(G)$ we have functorial identifications \beq\av_{B}^G\circ \HC(A)=\av_{B}^G(\eU\ast A)=\av_{B}^G(\eU)\ast A=\Spr\ast A\eeq
where $\Spr:=\av_{B}^G(\eU)=\Av_B^G(\eU)$ is the Springer sheaf. More generally for $C\in\Hecke,\ A\in \D_G(G)$ we have the functorial identifications
\beq\label{eq:avHC}
\av_{B}^G(C)\ast A \cong \av_{B}^G(C\ast \HC(A))\cong \av_{B}^G(C\ast A).\eeq

By the main results of \cite{BBM:04, BBM22}, for each $A\in \D_G(G)$ we have a natural isomorphism $\eL\ast A\cong \eL\ast_\ast A$, where $\ast_*$ denotes the convolution without compact support. Using this we can also obtain a left adjoint to $\HC_\L$ (where $d=\dim(G/U)$):
\beq
\av_U^G[2d](d):\eL\D(G)\eL\subset \D_U(G)\to \D_G(G)
\eeq
giving us the adjoint triple
\beq
(\av_U^G[2d](d),\HC_\L,\Av_U^G).
\eeq

\subsection{Main results}\label{sec:mr}
Let us now describe the main results of this paper. Recall that $N\subset G$ denotes the normalizer of the maximal torus $T$ and $W=N/T$ is the Weyl group relative to $T$. The group $W$ acts on the commutative algebraic group $T$ and hence acts on the symmetric monoidal category $\D(T)$ as well as on the braided monoidal category $\D_T(T)$ (see \cite[Appendix B]{BoDr:14}). Let us consider their $W$-equivariantizations, namely the categories $\D_W(T)$ and $(\D_T(T))^W=\D_N(T)$ of $W$-equivariant and $N$-equivariant complexes on the maximal torus $T$. It follows that the category $\D_W(T)$ is symmetric monoidal whereas $\D_N(T)$ is braided monoidal which is not symmetric in general. Note that if we let $N':=W\ltimes T$, we have the natural equivalences of braided monoidal categories: $\D_{N'}(T)\cong (\D_T(T))^W\cong \D_N(T)$; in fact, both categories are identified with the same full subcategory
in $\D_{\tilde N}(T)$ where $\tilde N = \tilde W\ltimes T$ for a finite group $\tilde W$
with a surjection $\tilde W \to W$ splitting the extension $T\to N \to W$. 

We refer to \S\ref{sec:centralsheaves} where the thick triangulated monoidal subcategories of central complexes on $T$, $\Dcent_W(T)\subset \D_W(T)$ and $\Dcent_N(T)\subset \D_N(T)$, are defined. These subcategories are defined in terms of actions of certain subgroups of $W$ on the cohomology with compact supports of Kummer twists of objects of $\D_W(T)$ or $\D_N(T)$. Equivalently, these may also be defined in terms of the Mellin transform. This notion agrees with that defined in \cite{Chen:1902} after twisting the $W$-equivariance structure by the sign character of $W$.  As we will see, this notion is closely related to objects in the Drinfeld center of the Hecke category $\eU\D_B(G)$ which are supported on $B\subset G$, which justifies the terminology ``central complexes'' on $T$.

%In \S\ref{sec:vanishingsheavesonT} the thick triangulated monoidal subcategories of vanishing sheaves on $T$, $\Dcirc_W(T)\subset \D_W(T)$ and $\Dcirc_N(T)\subset \D_N(T)$, are defined. 
We denote the parabolic induction functor with respect to $T\subset B$ simply as $\ind:\D_T(T)\to \D_G(G)$.
Note that the pair $(U,\Qlcl)$ is a Heisenberg admissible pair for the group $B$ in the sense of \cite{Desh:16} and we have a braided monoidal equivalence $\D_T(T)\cong \eU\D_B(B)$ defined by $A\mapsto e_U\ast A$. The parabolic induction functor can equivalently be described in terms of the averaging functor $\av_B^G=\Av_B^G$ as the composition \[\ind:\D_T(T)\cong\eU\D_B(B)\subset \eU\D_B(G)\xto{\Av_B^G}\D_G(G).\]
Hence by \cite[Lem. 7.1, Cor. B.47]{BoDr:14}, $\ind:\D_T(T)\to \D_G(G)$ has the structure of a weakly braided semigroupal functor between the two braided monoidal categories, namely we have functorial morphisms
\[
\ind(C)\ast \ind(D)\to \ind(C\ast D)\mbox{ for } C,D\in \D_T(T)\] 
which are compatible with the associativity and braiding constraints of $\D_T(T)$ and $\D_G(G)$.

We also have the parabolic restriction functor $\res:\D_G(G)\to \D_T(T)$ with respect to $T\subset B$ which gives us the adjoint pair of functors $(\res,\ind)$. Parabolic restriction can be described in terms of the Harish-Chandra restriction functor as the composition (where we let $|_B$ denote the restriction to $B\subset G$)
\beq\label{eq:resandHC}\res:\D_G(G)\xto{\HC}\eU\D_B(G)\xto{(\cdot)|_B}\eU\D_B(B)\cong \D_T(T),\eeq\beq A\mapsto (\eU\ast A)|_T[-2\dim U](-\dim U).\eeq
For each $A\in \D_G(G)$ we have the natural adjunction morphism
\beq\eU\ast A\cong\HC(A)\to \HC(A)_B\cong\eU\ast \res(A),\eeq
where $\HC(A)_B$ denotes the extension by zero of the restriction $\HC(A)|_B$. In \S\ref{sec:vanishingsheavesonG} we define the full subcategory $\Dcirc_G(G)\subset \D_G(G)$ of vanishing sheaves whose objects are those $A\in \D_G(G)$ whose Harish-Chandra transform $\HC(A)\in \eU\D_B(G)\eU$ is supported on the closed double coset $B\subset G$, i.e. those $A$ such that the adjunction morphism above is an isomorphism. 

In general at the level of triangulated categories there is no $W$-action on the parabolic restriction functor $\res:\D_G(G)\to \D_T(T)$ as explained in \cite{Gun:17}. However in \S\ref{sec:vanishingsheavesonG} we will see that when restricted to $\Dcirc_G(G)$, we can upgrade $\res$, which essentially agrees with  Harish-Chandra restriction, to a functor $$\res^\circ:\Dcirc_G(G)\to \D_N(T).$$ We now state our first main result:

\bthm\label{thm:main1}
The functor $\res^\circ$ above gives a triangulated braided monoidal equivalence
\[\res^\circ:\Dcirc_G(G)\rar{\cong} \Dcirc_N(T)\subset \D_N(T).\]
Moreover, this is a t-exact equivalence for the perverse t-structures on the source and target category.
\ethm

After defining all the relevant notions and establishing some preliminaries in \S\ref{sec:centralandvanishing}, we will prove the above result as Theorem \ref{thm:dcircg1} in \S\ref{sec:proofofthm:main1}. Note that since the parabolic restriction functor is t-exact by \cite{BYD:18}, the t-exactness part of the result follows.

To prove this result, we will use the vanishing conjecture stated in \cite{Chen:1902}, and proved in \cite{Chen:1909, BITV}. Since the adjoint action of $T$ on itself is trivial, we have a natural functor $\D(T)\to \D_T(T)$ which equips an object of $\D(T)$ with the trivial $T$-equivariance structure. In terms of $\Qlcl$-complexes on stacks, this is the pullback functor along the morphism $T/\ad T\cong T\times (\pt/T)\to T$. To describe this vanishing result, let us define the non-$T$-equivariant induction functor as the composition 
$$\ind:\D(T)\to \D_T(T)\to \D_G(G).$$ The $W$-action on $T$ gives rise to a symmetric monoidal action of $W$ on $\D(T)$ and the parabolic induction functor (from the  non-$T$-equivariant category $\D(T)$) naturally commutes with this $W$-action as described in \cite[Prop. 3.2(1)]{Chen:1902}. Hence the induction functor can be upgraded to a functor (in the non-$T$-equivariant setting)\[\ind:\D_W(T)\to \D^W_G(G),\]
where $\D^W_G(G)$ is the category of objects of $\D_G(G)$ equipped with a $W$-action. This is a natural generalization of Springer's $W$-action of the Springer sheaf $\Spr\in \D_G(G)$, which is by definition the parabolic induction of the unit $\delta_1\in \D_W(T)$. In \S\ref{sec:indresWaction} we will study the $W$-invariant summand of the induction functor 
\[\ind^W:\D_W(T)\to \D_G(G)\]
which is also a weakly braided semigroupal functor. We will prove that the $W$-invariant parabolic induction functor, when  restricted to the full subcategory of $\Dcirc_W(T)\subset \D_W(T)$ of central sheaves, defines a  braided monoidal functor $\ind^W:\Dcirc_W(T)\to \D_G(G)$. 
%Similarly we can define the full subcategory $\Dcirc_N(T)\subset \D_N(T)$ of vanishing sheaves.

The vanishing conjecture states that the $W$-invariant parabolic induction functor maps a central object in $\D_W(T)$ to a vanishing object in $\D_G(G)$, i.e. we have a braided monoidal functor
\[
\ind^W:\Dcirc_W(T)\to \Dcirc_G(G). 
\]

Next, we will proceed to study the bi-Whittaker category $\Wh$. Our next main result establishes an equivalence between the bi-Whittaker category and the category of $W$-equivariant central complexes on the torus $T$: 

\bthm\label{thm:main2}
There is a triangulated monoidal equivalence
\beq\label{eq:whittakertocentral}
\xi:\Wh\rar{\cong} \Dcirc_W(T)
\eeq
whose inverse is  given by the composition
\[\Dcirc_W(T)\xto{\ind^W}\Dcirc_G(G)\xto{\HC_\L}\Wh.\]
In particular the bi-Whittaker category $\Wh$ has the structure of a triangulated symmetric monoidal category. Moreover, the equivalence above is a t-exact equivalence for the perverse t-structure on $\Dcirc_W(T)$ and the perverse t-structure shifted by $\dim U$ on $\Wh$.
\ethm
The t-exactness in the statement above follows from the t-exactness of the functor $\HC_\L$ proved in \cite{BBM:04,BBM22} along with the t-exactness from Theorem \ref{thm:main1}. In order to construct the functor (\ref{eq:whittakertocentral}), in \S\ref{sec:bimodulecategory} we will equip the category $\D(T)$ with the structure of a $\eUop\D(G)\eUop-\Wh$-bimodule category such that the full monoidal subcategory $\D(T)\cong\eUop\D(\Bop)\subset \eUop\D(G)\eUop$ acts on the left on the bimodule category $\D(T)$ by the usual left convolution. Moreover we will define the objects $\Kcal_w\in \eUop\D(G)\eUop, w\in W$ whose left action on $\D(T)$ coincides with the adjoint action of $w$ on $\D(T)$. We will construct the functor $\xi$ and prove Theorem \ref{thm:main2} in \S\ref{sec:ptm2}.

In order to equip the category $\D(T)$ with the structure of a $\eUop\D(G)\eUop-\Wh$-bimodule category we will use the following easy lemma:
\blem\label{lem:bimodule}
The objects of the full subcategory $\eUop\D(G)\eL\subset \D(G)$ are supported on the open subset $\Uop T U\subset G$ and we have natural equivalences of triangulated categories 
\beq\label{eq:triangulatedequivalenceofbimodule}
\D(T)\cong \eUop\D(\Bop)\cong \D(B)\eL\cong \eUop\D(\Uop TU)\eL=\eUop\D(G)\eL.\eeq The category $\eUop\D(G)\eL$, and hence the category $\D(T)$ which is equivalent to it, has the structure of a $\eUop\D(G)\eUop-\Wh$-bimodule category given by convolution with compact supports on the left and right. The full monoidal subcategory $\D(T)\cong\eUop\D(\Bop)\subset \eUop\D(G)\eUop$ acts on the left on the bimodule category $\D(T)$ by the usual left convolution of $\D(T)$ on itself.
\elem
In \S\ref{sec:bimodulecategory} we prove this lemma and study the structure of the bimodule category $\D(T)$. In order prove Theorem \ref{thm:main2} we will study the left action of the Yokonuma-Hecke category $\eUop\D(G)\eUop$ on $\D(T)$.

\begin{rk}\label{Soe_rem}
The category of monodromic sheaves
in the Yokonuma-Hecke category $\D(\Uop\backslash G /\Uop)_{mon}$ was intensively studied, partly due to its relation to highest weight modules via the Localization Theorem of Beilinson-Bernstein and Brylinski-Kashiwara. A powerful method of its study pioneered by Soergel 
in \cite{Soe} (see also \cite{BeYu}) can be described from our
present perspective as follows. One considers the action of 
$\D(\Uop\backslash G /\Uop)$ on the category of monodromic sheaves $\eUop\D(G)\eL$, carrying a commuting right action of $\eL\D(G)\eL$. For concreteness, let us restrict attention  to monodromic sheaves with unipotent monodromy; then Theorem \ref{thm:main2} and Lemma \ref{lem:bimodule} show that the corresponding subcategory in $\eL\D(G)\eL$ acting on $\eUop\D(G)\eL$
is identified with $D^b(Coh_0( {\mathfrak{t}}/W))$ acting on $D^b(Coh_0({\mathfrak{t}}))$, here $Coh_0$ denotes the category of coherent
sheaves set-theoretically supported at 0. The commuting action 
of the monodromic Yokonuma-Hecke category leads to a monoidal functor
from the latter to $D^b(Coh({\mathfrak{t}}\times_{{\mathfrak{t}}/W} {\mathfrak{t}}) $, coherent sheaves on Soergel variety known also as Soergel
bimodules.
Notice that in \cite{BeVi} Soergel's method was generalized to the category of $({\mathfrak{g}},K)$-modules where $K$ is the complexification of a maximal
compact subgroup in a quasi-split real from of a complex reductive group $G$.
The functor to coherent sheaves on the so called block variety playing a
central role in that work can also 
be described in terms of Whittaker sheaves.

It would be interesting to use our results to extend Soergel's method to
a description of the full (rather than monodromic) Yokonuma-Hecke category.
In the $D$-modules setting this is addressed in Gannon's work \cite{Gan23}.
\end{rk}

In \S\ref{sec:yokonuma} we study the Yokonuma-Hecke category and certain special objects in it defined by Kazhdan and Laumon in \cite{KazLau}. Namely, for each $w\in W$ we have a shifted perverse sheaf $\Kcal_w\in \eUop\D(G)\eUop$ whose left action on $\D(T)$ coincides with the adjoint action of $w$ on $\D(T)$:
\bthm\label{thm:jsheaves} (\cite{KazLau}.)
For each $w\in W$ we have a distinguished object $\Kcal_w\in \eUop\D(G)\eUop$, with $\Kcal_{1_W}=\eUop$, such that we have the following:\\
(i) The object $\Kcal_w[-\dim \Uop]$ is a simple perverse sheaf. It is the middle extension of its restriction to $\Bop w\Bop$ where it is an irreducible local system.\\
(ii) If $w_1,w_2\in W$ are such that $l(w_1w_2)=l(w_1)+l(w_2)$ then we have a natural isomorphism $$\Kcal_{w_1w_2}\cong \Kcal_{w_1}\ast \Kcal_{w_2}.$$
(iii) The left action of $\Kcal_w$ on the bimodule category $\D(T)$ coincides with the adjoint action of $w$ on $\D(T)$.
\ethm
\noindent The statements (i) and (ii) are proved in \cite{KazLau}. Statement (iii) will be proved in \S\ref{sec:bimodulecategory}.

As an immediate consequence of our main results we obtain:
\bcor\label{equiv_ab}
Theorems \ref{thm:main1} and \ref{thm:main2} give us equivalences of the abelian categories
\[\Perv^\circ_G(G)\cong \Perv^\circ_W(T)\cong \Perv(\Wh) %[\dim U]
,\]
where $\Perv(\Wh):=\Wh\cap \Perv(G)$. 
\ecor

\begin{rk}
As follows from the above, the equivalence $\Perv^\circ_G(G)\to \Perv(\Wh)$
is given by the natural averaging functor (up to a homological shift by $\dim U$).
However, the natural averaging functor $\Perv(\Wh)\to \D_G(G)$ does not yield the 
inverse equivalence. For unipotent sheaves, the latter functor is studied in \cite{BeTo}
where it is related to Hochschild homology for Soergel bimodules.
\end{rk}

\subsection*{Acknowledgments}
We thank Tsao-Hsien Chen, Tom Gannon, Sam Gunningham, Kostiantyn Tolmachov and Alexander Yom Din for helpful correspondence. RB was partly supported by NSF grant DMS-2101507. TD was supported by the Department of Atomic
Energy, Government of India, under project no.12-R\&D-TFR-5.01-0500.

\section{Central sheaves and vanishing sheaves}\label{sec:centralandvanishing}
In this section we introduce the categories of central sheaves on the torus $T$ and prove some results about them. We also introduce the category of vanishing sheaves on $G$ and prove that it is equivalent to the category of ($N$-equivariant) central sheaves on $T$.

\subsection{Central sheaves on the torus}\label{sec:centralsheaves}
In this section we define certain full subcategories of the categories $\Dcent_W(T)\subset \D_W(T)$ and $\Dcent_N(T)\subset \D_N(T)$ which we call the categories of central sheaves on the torus. Our definition differs from the one given in \cite{Chen:1902,Chen:1909} by the sign character.

 Following \cite{GL}, let $\Ccal(T)$ denote the $\Qlcl$-scheme whose $\Qlcl$-points parametrize the multiplicative (or Kummer) $\Qlcl$-local systems on $T$, i.e. continuous $\Qlcl$-valued characters of the tame \'etale fundamental group $\pi_1^t(T)$. Note that we have a canonical identification $\pi_1^t(T)=X_*(T)\otimes_{\bbZ}\widehat{\bbZ}(1)$. 
 We will use the same notation, $\Ccal(T)$, to denote the set of $\Qlcl$-points of the above scheme. 
 
 The action of $W$ on $T$ induces an action of $W$ on $\Ccal(T)$ defined by $w(\chi)={w^{-1}}^*\chi$ for $w\in W,\chi\in \Ccal(T)$. For $\chi\in \Ccal(T)$, let $W_\chi\subset W$ denote the stabilizer of $\chi$. Note that if $w\in W_\chi$ then we in fact have a canonical isomorphism of multiplicative local systems $w(\chi)\cong \chi$. If $\alpha^\vee:\Gm\to T$ is a coroot such that ${\alpha^\vee}^*\chi\cong \Qlcl$, then $s_\alpha(\chi)\cong \chi$, i.e. $s_\alpha\in W_\chi$, where $s_\alpha\in W$ is the reflection corresponding to the root $\alpha$. Consider the following subset of $\Phi\times \Ccal(T)$
 \beq
 \Pcal:=\{(\alpha,\chi)|{\alpha^\vee}^*\chi\cong\Qlcl\}.
 \eeq

\bdefn\label{def:Omega}
For $\chi\in \Ccal(T)$, let $W^\circ_\chi\subset W_\chi$ be the subgroup generated by the $s_\alpha$ such that $(\alpha,\chi)\in \Pcal$. Then $W^\circ_\chi$ is a normal subgroup of $W_\chi$ (see Lemma \ref{lem:Omega} below). We let $\Omega_\chi$ denote the quotient group $W_\chi/W^\circ_\chi$. It is an abelian group by Lemma \ref{lem:Omega} below.
\edefn 
Note that for $A\in \D_W(T)$, $\chi\in \Ccal(T)$ and $w\in W$ we have functorial isomorphisms \beq
\phi_{w}^{A,\chi}:H^*_c(T,A\otimes \chi)\rar{\cong} H^*_c(T,A\otimes w(\chi))
\eeq satisfying the following cocycle condition for $w_1,w_2\in W$:
\beq\label{eq:cocyclecondition}
\phi_{w_1w_2}^{A,\chi}=\phi_{w_1}^{A,w_2(\chi)}\circ \phi_{w_2}^{A,\chi}.
\eeq
Hence if $w\in W_\chi$ we have the automorphism
\beq\label{eq:actionofstabilizer}
\phi_{w}^{A,\chi}:H^*_c(T,A\otimes \chi)\rar{\cong} H^*_c(T,A\otimes \chi)
\eeq
which gives rise to an action of $W_\chi$ on $H^*_c(A\otimes\chi)$.
In particular for $(\alpha,\chi)\in \Pcal$ we have the automorphism $\phi_{s_\alpha}^{A,\chi}:H^*_c(T,A\otimes \chi)\rar{\cong} H^*_c(T,A\otimes \chi).$
\blem\label{lem:centralconditions}
Let $A\in \D_W(T)$. Then the following conditions are equivalent:\\
(i) For each $\chi\in \Ccal(T)$, the subgroup $W^\circ_\chi$ of $W_\chi$ acts trivially on $H^*_c(T,A\otimes\chi)$ via Equation (\ref{eq:actionofstabilizer}).\\
(ii) For each $(\alpha,\chi)\in \Pcal$, the automorphism $\phi_{s_\alpha}^{A,\chi}$  of $H^*_c(T,A\otimes\chi)$ is the identity.\\
(iii) For each $(\alpha,\chi)\in \Pcal$ such that $\alpha\in \Delta$ is a simple root, the automorphism $\phi_{s_\alpha}^{A,\chi}$  of $H^*_c(T,A\otimes\chi)$ is the identity.\\
(iv) For each simple root $\alpha\in \Delta$, the forgetful image of the object $A$ in $\D_{\langle s_{\alpha}\rangle}(T)$ lies in $\Dcirc_{\langle s_{\alpha}\rangle}(T)$.
\elem
\bpf
It is clear that (i) $\iff$ (ii) $\implies$ (iii) $\iff$ (iv), so let us prove that (iii)$\implies$(ii). Note that we have the action of $W$ on $\Pcal$ given by $w\cdot(\alpha,\chi)=(w(\alpha),w(\chi))$. Also we have $s_{w(\alpha)}=ws_\alpha w^{-1}$. Hence using Equation (\ref{eq:cocyclecondition}) for each $(\alpha,\chi)\in \Pcal, w\in W$ we have the commutative diagram
\beq
\xymatrixcolsep{3pc}\xymatrix{
H^*_c(T,A\otimes\chi)\ar[d]_{\phi^{A,\chi}_w}\ar[r]^-{\phi^{A,\chi}_{s_\alpha}} & H^*_c(T,A\otimes\chi)\ar[d]^{\phi^{A,\chi}_w}\\
H^*_c(T,A\otimes w(\chi))\ar[r]^-{\phi^{A,w(\chi)}_{s_{w(\alpha)}}} & H^*_c(T,A\otimes w(\chi)).
 }
\eeq
Hence we see that condition (ii) is satisfied for the pair $(\alpha,\chi)\in \Pcal$ if and only if it is satisfied for the pair $(w(\alpha),w(\chi))$. Now the implication (iii)$\implies$(ii) follows, since for any $\alpha\in \Phi$ we can find  $w\in W$ such that $w(\alpha)\in \Delta$. 
\epf

Note that the forgetful functor $\D_T(T)\to \D(T)$ is monoidal and $W$-equivariant for the natural $W$-actions on the categories $\D_T(T), \D(T)$. Hence after $W$-equivariantization, we have a natural monoidal forget functor $\D_N(T)=\D_T(T)^W\to \D_W(T)=\D(T)^W$.
\bdefn\label{def:central}
An object $A\in \D_W(T)$ is said to be central if it satisfies any of the equivalent conditions of Lemma \ref{lem:centralconditions}. The full subcategory of central objects in $\D_W(T)$ is denoted by $\Dcent_W(T)$. Similarly, we say that an object $A\in \D_N(T)$ is central if its forgetful image in $\D_W(T)$ is central. The full subcategory of central objects in $\D_N(T)$ is denoted by $\Dcent_N(T)$.
\edefn

\blem\label{lem:centralthickmon}
The full subcategories $\Dcirc_W(T)\subset \D_W(T)$ and $\Dcirc_N(T)\subset \D_N(T)$ are thick triangulated monoidal subcategories.
\elem
\bpf
Since the forgetful functor $\D_N(T)\to \D_W(T)$ is triangulated and monoidal, it is enough to prove the lemma for the full subcategory $\Dcirc_W(T)\subset \D_W(T)$. By definition it is clear that $\Dcirc_W(T)$ is closed under direct summands. Now suppose that $A_1\to A_2\to A_3\to$ is a distinguished triangle in $\D_W(T)$ such that $A_1,A_2\in \Dcirc_W(T)$. Then for each $\chi\in \Ccal(T)$ we have the distinguished triangle $$H^*_c(T,A_1\otimes\chi)\to H^*_c(T,A_2\otimes\chi)\to H^*_c(T,A_3\otimes\chi)\to$$ in $D^b(\Rep(W_\chi))$, where we know that the $W_\chi^\circ$-actions on the first two cohomologies are trivial. Hence the $W_\chi^\circ$-action on $H^*_c(T,A_3\otimes \chi)$ must also be trivial, i.e. $A_3\in \Dcirc_W(T)$. Hence $\Dcirc_W(T)\subset \D_W(T)$ is a thick triangulated subcategory. 

To prove that it is a monoidal subcategory, first note that the unit $\delta_1\in \D_W(T)$ lies in $\Dcirc_W(T)$. Now let $A_1,A_2\in \Dcirc_W(T)$ and $\chi\in \Ccal(T)$. In particular, $\chi$ is a multiplicative local system on $T$, namely $\mu^*\chi\cong \chi\boxtimes\chi$, where $\mu:T\times T\to T$ is the multiplication map. Now we have 
\[H^*_c(T,(A_1\ast A_2)\otimes \chi)=H^*_c(T,\mu_!(A_1\boxtimes A_2)\otimes \chi)=H^*_c(T,\mu_!((A_1\boxtimes A_2)\otimes \mu^*\chi))\]
\[=H^*_c\left(T,\mu_!((A_1\otimes\chi)\boxtimes (A_2\otimes \chi))\right)=H^*_c\left(T\times T,(A_1\otimes\chi)\boxtimes (A_2\otimes \chi)\right)\]
\[=H^*_c(T,A_1\otimes\chi)\otimes H^*_c(T,A_2\otimes\chi).\]
Then since the $W^\circ_\chi$-actions on $H^*_c(T,A_1\otimes\chi)$ and $H^*_c(T,A_2\otimes\chi)$ are trivial, it follows that the $W^\circ_\chi$-action on $H^*_c(T,(A_1\ast A_2)\otimes\chi)$ is trivial. Hence $A_1\ast A_2\in \Dcirc_W(T)$, proving that it is a monoidal subcategory.
%That it is a monoidal subcategory follows from \cite[Cor 4.3]{Chen:1902} and the simple observation that the unit $\delta_1\in \D_W(T)$ lies in $\Dcirc_W(T)$.
\epf

\blem\label{lem:tstructurecentral}
The restrictions of the perverse t-structures on the categories $\D_N(T)$ and $\D_W(T)$ define t-structures on the categories $\Dcirc_N(T)$ and $\Dcirc_W(T)$ of central sheaves on the torus.
\elem
\bpf
This follows from the proof of Lemma 7.2 in \cite{Chen:1909} which proves that the perverse truncations of central sheaves are also central.
\epf

\brk\label{rk:pervtstructure}
In other words we have well-defined perverse t-structures on the full subcategories $\Dcirc_N(T)$ and $\Dcirc_W(T)$. We denote their hearts by $\Perv^\circ_N(T):=\Perv_N(T)\cap \Dcirc_N(T)$ and $\Perv^\circ_W(T):=\Perv_W(T)\cap \Dcirc_W(T)$ respectively. We have equivalences of hearts $\Perv_N(T)\cong \Perv_W(T)$, and hence $\Perv^\circ_N(T)\cong \Perv^\circ_W(T)$.
\erk

\subsection{Monodromic sheaves and Mellin transform}\label{sec:monodromicmellin}
Recall that we have an action of the Weyl group $W$ on the scheme $\Ccal(T)$ parametrizing the $\Qlcl$-multiplicative local systems on $T$. Also recall that for a multiplicative local system $\chi\in \Ccal(T)$ we have its stabilizer $W_\chi$ and the normal subgroup $W^\circ_\chi\leq W_\chi$. Moreover, we have a semidirect product decomposition (see \cite[\S1.11]{Lus:NUCS}) $W_\chi=W^\circ_\chi\rtimes\Omega_\chi$ with the quotient $W_\chi/W^\circ_\chi=\Omega_\chi\hookrightarrow W_\chi$ being abelian (see also Lemma \ref{lem:Omega} below). In particular, if $\chi,\chi'\in \Ccal(T)$ lie in the same $W$-orbit $\Theta$ then we have a canonical identification of the abelian quotients $\Omega_\chi=\Omega_{\chi'}$ which we will denote as $\Omega_{\Theta}$. 

We also recall from \cite[\S2]{Desh:17} that for each $\chi\in \Ccal(T)$ we have the minimal quasi-idempotent $e_\chi:=\chi\otimes\omega_T=\chi[2\dim T](\dim T)\in \D(T)$ and the corresponding thick triangulated ideal $\D(T)_\chi\subset \D(T)$ of $\chi$-monodromic complexes on $T$. Moreover by \cite{GL}, we have the Mellin transform (with compact supports) 
$$\Mcal_!:(\D(T),\ast)\rar{} (D^b\Coh\Ccal(T),\otimes)$$ 
which is a monoidal functor and which induces an equivalence 
$$\Mcal_!:\D(T)_\chi\rar{\cong} D^b\Coh_{\chi^{-1}}\Ccal(T)$$
where $D^b\Coh_{\chi^{-1}}\Ccal(T)\subset  D^b\Coh\Ccal(T)$ is the full subcategory of complexes of coherent sheaves set theoretically supported on the closed point ${\chi^{-1}}\in \Ccal(T)$.
% Now suppose that $D\subset \Ccal(T)$ be a nilpotent thickening of $\{\chi\}$, i.e. a closed subscheme of $\Ccal(T)$ with $D_{red}=\{\chi\}$. Note that we have a natural functor $D^b\Coh D\to D^b\Coh_\chi\Ccal(T)$. Let $\Ecal_D\in \D(T)_\chi$ be the object such that $\Mcal_!(\Ecal_D)=\Ocal_D$, the structure sheaf of $D\subset \Ccal(T)$.

The projections $\Ccal(T)\to \Ccal(T)//W_\chi^\circ\to \Ccal(T)//W_\chi$ give us the pullback functors $D^b\Coh(\Ccal(T)//W_\chi)\to D^b\Coh(\Ccal(T)//W^\circ_\chi)\to D^b\Coh\Ccal(T)$ which are symmetric monoidal for the tensor product.
Let $\mathbbm{O}_{\chi^{-1}}^\circ$ (resp. $\mathbbm{O}_{\chi^{-1}}$) in $D^b\Coh_\chi\Ccal(T)\subset D^b\Coh\Ccal(T)$ be the pullback of the structure sheaf of the image of the point $\chi^{-1}$ in $\Ccal(T)//W_\chi^\circ$ (resp. in $\Ccal(T)//W_\chi$) under the above projections. We denote the structure sheaf of the point $\chi^{-1}\in \Ccal(T)$ as $\Ocal_{\chi^{-1}}$. Note that the Mellin transform $\Mcal_!$ of the minimal quasi-idempotent $e_\chi\in \D(T)$ is exactly the structure sheaf $\Ocal_{\chi^{-1}}$. Similarly let $\ecircchi$ (resp. $\echi$) in $\D(T)_\chi$ denote the inverse Mellin transform of $\mathbbm{O}_{\chi^{-1}}^\circ$ (resp. $\mathbbm{O}_{\chi^{-1}}$) in $D^b\Coh_\chi\Ccal(T).$ 

The objects $e_\chi,\ecircchi,\echi$ have natural $W_\chi$-equivariant structures. Also for any $C\in \D_{W_\chi}(T)$ we have a $W_\chi$ action on $H^*_c(T,C\otimes \chi^{-1})=1^*(C\ast\chi)$, where we are considering $\chi,\chi^{-1}$ as a multiplicative local system in $\D_{W_\chi}(T)$ and $1:\pt\to T$ is the unit of $T$. Consider the object $1_*H^*_c(T,C\otimes \chi^{-1})\in \D_{W_\chi}(1)\subset \D_{W_\chi}(T)$. By using the properties of the Mellin transform and the results from \cite{Ne:08} about descent of coherent sheaves we obtain (see \cite[Prop. 5.4 and its proof]{Chen:1909}):
\bprop\label{p:convwithechi}
Let $C\in \D_{W_\chi}(T)$. \\
(i) We have an isomorphism $C\ast e_\chi\cong 1_*H^*_c(T,C\otimes \chi^{-1})\ast e_\chi$ in $\D_{W_\chi}(T)$.\\
(ii) If the action of the subgroup $W^\circ_{\chi}$ on $H^*_c(T,C\otimes\chi^{-1})$ is trivial, then we have an isomorphism 
$$C\ast \ecircchi \cong 1_*H^*_c(T,C\otimes \chi^{-1})\ast \ecircchi \mbox{ in }\D_{W_\chi}(T).$$
In particular we have $\ecircchi\ast \ecircchi \cong 1_*H^*_c(T,\ecircchi\otimes \chi^{-1})\ast \ecircchi \mbox{ in }\D_{W_\chi}(T).$\\
(iii) If the action of $W_{\chi}$ on $H^*_c(T,C\otimes\chi^{-1})$ is trivial, then we have an isomorphism 
$$C\ast \echi \cong 1_*H^*_c(T,C\otimes \chi^{-1})\ast \echi \mbox{ in }\D_{W_\chi}(T).$$
In particular we have $\echi\ast \echi \cong 1_*H^*_c(T,\echi\otimes \chi^{-1})\ast \echi \mbox{ in }\D_{W_\chi}(T).$
\eprop

According to our definition of central sheaves on $T$, we are most interested in those $C\in \D_{W_\chi}(T)$ such that the $W^\circ_\chi$-action on $H^*_c(T,C\otimes\chi^{-1})$ is trivial. For such a $C$, we have an action of $\Omega_\chi=W_\chi/W^\circ_\chi$ on $H^*_c(T,C\otimes\chi^{-1})$ and taking the isotypic components we obtain a canonical decomposition 
\beq
1_*H^*_c(T,C\otimes\chi^{-1})=\bigoplus\limits_{\rho\in \widehat{\Omega}_{\chi}}H^*_c(T,C\otimes\chi^{-1},\rho)\otimes \delta_1^\rho \mbox{ in } \D_{W_\chi}(T)
\eeq
where for any irreducible character $\rho\in\widehat{\Omega}_\chi$, we define $H^*_c(T,C\otimes\chi^{-1},\rho):=\RHom_{\Omega_\chi}(\rho,H^*_c(T,C\otimes\chi^{-1}))\in D^b\Vec_{\Qlcl}$ and we define $\delta^\rho_1\in \D_{W_\chi}(1)\subset \D_{W_\chi}(T)$ to be the delta sheaf supported at $1\in T$ equipped with the $W_\chi$-equivariant structure corresponding to $\rho:W_\chi\to\Omega_\chi\to \Qlcl^\times$.

We now set ${\ecircchi}^\rho:=\delta_1^\rho\ast \ecircchi\in \D_{W_\chi}(T)$, i.e. it has the same underlying object in $\D(T)$ as $\ecircchi$, but with the $W_\chi$-equivariance structure defined  by the character $\rho$. By Proposition \ref{p:convwithechi}(ii) if $C\in \D_{W_\chi}(T)$ is such that the $W^\circ_\chi$-action on $H^*_c(T,C\otimes\chi^{-1})$ is trivial then we have
\beq
C\ast \ecircchi \cong \bigoplus\limits_{\rho\in\widehat{\Omega}_\chi}H^*_c(T,C\otimes \chi^{-1},\rho)\otimes {\ecircchi}^\rho\mbox{ in }\D_{W_\chi}(T).
\eeq
If $\Theta\subset \Ccal(T)$ denotes the $W$-orbit of $\chi$, we have the equivalence $D^b_W\Coh(\Theta)\cong D^b\Rep W_\chi$. Hence if $C\in \Dcirc_W(T)$ is a central sheaf on $T$ and  if $\chi,\chi' \in\Ccal(T)$ are in the same $W$-orbit then we have a canonical identification $\RHom_{\Omega_\chi}(\rho,H^*_c(T,C\otimes\chi^{-1}))\cong \RHom_{\Omega_{\chi'}}(\rho',H^*_c(T,C\otimes{\chi'}^{-1}))$, where we recall that we have a canonical identification $\Omega_\chi\cong\Omega_{\chi'}$ and $\rho'$ corresponds to $\rho$ under this identification. Let us denote this previous complex of vector spaces as $H^*_c(T,C,\Theta^{-1},\rho)$ for any $C\in \Dcirc_W(T)$, $W$-orbit $\Theta\subset \Ccal(T)$, $\rho\in \widehat{\Omega}_\Theta$.
\bcor\label{cor:convwithecirctheta}
Let $\Theta\subset \Ccal(T)$ be the $W$-orbit of $\chi\in \Ccal(T)$ and let $\rho:\Omega_\Theta\to \Qlcl^\times$ be a character. Define ${\bbecirc_\Theta}^\rho:=\av_{W_\chi}^W{\ecircchi}^\rho\in \Dcirc_W(T)$. Then for any $C\in \Dcirc_W(T)$ we have
$$C\ast \bbecirc_\Theta \cong \bigoplus\limits_{\rho\in\widehat{\Omega}_\chi}H^*_c(T,C, \Theta^{-1},\rho)\otimes {\bbecirc_\Theta}^\rho\mbox{ in }\Dcirc_{W}(T).$$
In particular $\bbecirc_\Theta\ast \bbecirc_\Theta \cong \bigoplus\limits_{\rho\in\widehat{\Omega}_\chi}H^*_c(T,\bbecirc_\Theta, \Theta^{-1},\rho)\otimes {\bbecirc_\Theta}^\rho\mbox{ in }\Dcirc_{W}(T).$
\ecor

\subsection{Comparison with groups having connected center}
For a reductive group $G$ let $\bar{G}:=G^{\opname{ab}}\times G^{\ad}=G/[G,G]\times G/Z$, where $Z\leq G$ is the center of $G$ and $[G,G]\leq G$ is the commutator subgroup. Note that we have $G=Z\cdot[G,G]$ and that $\Lambda:=Z\cap[G,G]$ is finite. Hence we have the central extension
\beq\label{eq:centext}
1\to \Lambda \to G \xto{\varpi_G} \bar{G}\to 1.
\eeq
Also note that $G^{\ab}=G/[G,G]=Z/\Lambda$ is a torus. If $T$ is a maximal torus of $G$, then we let $T^{\ad}:=T/Z$ denote the corresponding maximal torus of $G^{\ad}$, and $\bar{T}:=T/\Lambda=G^{\ab}\times T^{\ad}$ the corresponding maximal torus of $\bar{G}$. The normalizer of $\bar{T}$ is $\bar{N}=N/\Lambda$ and we have a canonical identification of the Weyl groups of $G$ and $\bar{G}$. We obtain the $W$-equivariant central isogeny (where $W$ acts trivially on $\Lambda$):
\beq
1\to \Lambda\to T\xto{\varpi_T} \bar{T}\to 1.
\eeq
We let $\Lambda^*:=\Hom(\Lambda,\Qlcl^\times)$ be the Pontryagin dual of $\Lambda$. For each $\zeta\in \Lambda^*$ we can define a multiplicative local system $\chi_\zeta$ on $\bar{G}$ using the central extension (\ref{eq:centext}). We denote its restriction to $\bar{T}$ by $\chiTzeta\in \Ccal(\bar{T})$. This gives us an embedding $\Lambda^*\subset \Ccal(\bar{T})$.

Since $\varpi_T$ is an isogeny, the $W$-equivariant pullback morphism $\varpi_T^*:\Ccal(\bar{T})\to \Ccal(T)$ is surjective and is in fact a $\Lambda^*$-torsor. In other words given any multiplicative local system $\Qcal$ on $T$, the set $\Ccal(\bar{T})_{\Qcal}:=\{\chi\in\Ccal(\bar{T})|\varpi_T^*\chi\cong \Qcal\}$ is a $\Lambda^*$-torsor under tensor product of multiplicative local systems. The stabilizer of $\Qcal$ for the $W$-action, $W_{\Qcal}\leq W$ acts on the fiber $\Ccal(\bar{T})_{\Qcal}$ and for $\chi\in \Ccal(\bar{T})_{\Qcal}$ we have $W_\chi\leq W_\Qcal$. Then the map
\beq
W_{\Qcal}\ni w\mapsto w(\chi)\otimes\chi^{-1}\in \Lambda^*\hookrightarrow \Ccal(\bar{T})
\eeq
is independent of $\chi\in\Ccal({\bar{T}})_\Qcal$ and in fact defines a group homomorphism $W_\Qcal\to \Lambda^*$ whose kernel is exactly $W_\chi\leq W_\Qcal$. In particular the subgroup $W_\chi$ is itself independent of $\chi\in \Ccal({\bar{T}})_\Qcal$. In fact we have
\blem\label{lem:Omega}
As in the setting above, let $\Qcal\in\Ccal(T)$ and $\chi\in\Ccal(\bar{T})_\Qcal$. Then we have $W_\chi=W^\circ_{\chi}=W^\circ_\Qcal\leq W_\Qcal$. Hence for every $\Qcal\in \Ccal(T)$ the subgroup $W^\circ_\Qcal\leq W_\Qcal$ is normal and we have an inclusion 
\beq\label{eq:Omega}
\iota_\Qcal:W_\Qcal/W^\circ_{\Qcal}=\Omega_\Qcal\hookrightarrow \Lambda^*.
\eeq
\elem
\bpf
The center of $\bar{G}$ is the torus $G^{\ab}$, and hence is connected. Hence for any $\chi\in\Ccal(\bar{T})$ we have $W_\chi=W^\circ_\chi$ (see \cite[Thm. 5.13]{DL:76}). Now if $\alpha^\vee:\Gm\to T$ is a coroot of $G$, then the corresponding coroot of $\bar{G}$ is the composition $\bar{\alpha}^\vee:\Gm\xto{\alpha^\vee} T \xto{\varpi_T}\bar{T}$. Hence if $\chi\in\Ccal(\bar{T})_\Qcal$ then ${\alpha^\vee}^*\Qcal\cong \Qlcl$ if and only if ${\bar{\alpha}^\vee}^*\chi\cong \Qlcl$. Hence by Definition \ref{def:Omega} both $W^\circ_\Qcal$ and $W^\circ_\chi$ are generated by the same reflections in $W$, and hence are equal. This combined with the argument preceding the statement of the lemma gives us (\ref{eq:Omega}).
\epf

\subsection{Non-$T$-equivariant parabolic induction and restriction functors}\label{sec:indresWaction}
In this subsection we consider non-$T$-equivariant versions of the parabolic induction and restriction functors. Recall first that we have the parabolic induction functor 
$\ind=\ind_{T,B}^G:\D_T(T)\to \D_G(G)$ which can be identified with the following averaging functor restricted to the full subcategory $e_U\D_B(B)\subset \eU\D_B(G)=e_U\D_B(G)e_U$:
$$\av_B^G=\Av_B^G: e_U\D_B(G)\to \D_G(G).$$ Hence by \cite{BoDr:14} $\ind:\D_T(T)\cong \eU\D_B(B)\to \D_G(G)$ is a weakly braided semigroupal functor. Note that we have the equivalences $\D_T(T)\cong \eU\D_B(B)$ and $\D(T)\cong \eU\D(B)$ defined by $A\mapsto \eU\ast A$, and in the reverse direction by $C\mapsto C|_T[-2\dim U](-\dim U)$.

By a slight abuse of notation, we also denote by $\ind$ the non-$T$-equivariant parabolic induction functor, defined as the composition $$\ind:\D(T)\to\D_T(T)\to \D_G(G).$$ It also has a weakly braided semigroupal structure. We will also use the same notation to denote the functor $\ind:\D(T)\to \D(G)$ studied in \cite{Chen:1902}, and which is obtained by forgetting the $G$-equivariance structure. 

Consider the non-equivariant parabolic restriction functor (also denoted by $\res$) $$\res:\D(G)\to\eU\D(B)\cong \D(T),\  A\mapsto (\eU\ast A)|_T[-2\dim U](-\dim U).$$ 
\brk
Note that unlike in the equivariant setting, the non-equivariant parabolic restriction functor above is \emph{not adjoint} to the non-equivariant parabolic induction functor.
\erk
We have the commutative diagram of functors
\[
\xymatrix{
\D_G(G)\ar[r]^-{\res}\ar[d] & \D_T(T)\ar[d] \\
\D(G) \ar[r]^-{\res} & \D(T).
}
\]

By \cite{BYD:18} the parabolic induction and restriction functors are t-exact with respect to the perverse t-structure. Also note that we have a natural identifications of the hearts $\Perv(T)\cong \Perv_T(T)$ sitting inside $\D(T)$ and $\D_T(T)$ respectively, although these latter triangulated categories are quite different. Similarly we can consider $\Perv_G(G)$  as a full subcategory of $\Perv(G)$. We have an adjoint pair of exact functors
\[\res:\Perv_G(G)\to \Perv(T)\cong\Perv_T(T)\mbox{ and }\]\[\ind:\Perv_T(T)\cong \Perv(T)\to \Perv_G(G).\]

By \cite[\S3.2]{Chen:1902}, for each $w\in W$ we have a natural equivalence of functors
\beq\label{eq:indcircw}
\ind\circ w\cong \ind:\Perv(T)\to \Perv_G(G).
\eeq This defines an action of $W$ on the $W$-equivariant induction functor $\ind:\Perv_W(T)\to \Perv_G(G)$ from the category of  $W$-equivariant perverse sheaves, except that it will be convenient for us to consider the $W$-action which differs from {\it loc. cit.} by the sign representation of $W$. We have a canonical decomposition of the above functor as
\beq
\ind=\bigoplus\limits_{(\rho,V_\rho)\in\Irr W}V_\rho\otimes\ind^{W,\rho}:\Perv_W(T)\to \Perv_G(G).
\eeq
\bdefn\label{def:indW}
We denote the direct summand $\ind^{W,\tr}$ corresponding to taking $W$-invariants of the $W$-action by 
$$\ind^W:\Perv_W(T)\to \Perv_G(G).$$  We extend it to define a functor
$$\ind^W:\D_W(T)\to \D_G(G)$$
and call it the $W$-invariant parabolic induction functor. It has a weak braided semigroupal structure defined using the weak braided semigroupal structure on $\ind:\D_T(T)\to \D_G(G)$: $$\ind^W(A_1)\ast\ind^W(A_2)\to\ind(A_1)*\ind(A_2)\to \ind(A_1\ast A_2)\to \ind^W(A_1\ast A_2) \mbox{ for any }A_1,A_2\in \D_W(T).$$
\edefn
Consider the unit object (with trivial $W$-equivariant structure) $\delta_1\in \D_W(T)$. Then $\ind(\delta_1)=\Spr$ is the Springer sheaf in $\D_G(G)$. By our convention about the $W$-action on $\ind$ we get that \beq\ind^W\delta_1=\delta_1\in \D_G(G).\eeq

Note that using (\ref{eq:indcircw}) and adjunction, we get natural equivalences of functors $w\circ \res \cong \res:\Perv_G(G)\to \Perv(T)$ for each $w\in W$. 
\bdefn\label{defn:resw}
Using the above observation we get a functor $\res_W:\Perv_G(G)\to \Perv_W(T)$.
\edefn

\blem\label{lem:resWindW}
(i) The pair of functors $\left(\res_W:\Perv_G(G)\to \Perv_W(T), \ind^W:\Perv_W(T)\to \Perv_G(G)\right)$ is an adjoint pair.\\
(ii) For any $C\in \Perv_W(T)$ the adjunction morphism is an isomorphism
$$\res_W\circ\ind^W(C)\rar{\cong} C.$$
\elem
\bpf
Statement (i) is a formal consequence of the fact that $(\res,\ind)$ is an adjoint pair of functors and the natural equivalences from (\ref{eq:indcircw}). Statement (ii) is \cite[Prop. 3.4]{Chen:1902}, \cite[Thm. 2.7]{BK:03}.
\epf

\subsection{Vanishing sheaves on $G$}\label{sec:vanishingsheavesonG}
We will define a certain full subcategory $\Dcirc_G(G)\subset \D_G(G)$ of ``vanishing sheaves''. Recall that the Harish-Chandra transform is a central monoidal functor, namely it factors through the Drinfeld center as
\[
\D_G(G)\to \Zcal(e_U\D_B(G)e_U)\to e_U\D_B(G)e_U.
\]

The pair $(U,\Qlcl)$ is a Heisenberg admissible pair for $B$ in the sense of \cite{Desh:16} and we have an identification $\D_T(T)\cong e_U\D_B(B) \subset \eU\D_B(G)= e_U\D_B(G)e_U.$ 
Explicitly, the identification is given by
\beq\label{eq:identificationBT}
\D_T(T)\ni A \mapsto e_U\ast A\in e_U\D_B(B) \mbox{ and }\eU\D_B(B)\ni C\mapsto C|_{T}[-2\dim U](-\dim U)\in \D_T(T).
\eeq

Note that for $A\in \D_G(G)$ we have the natural adjunction morphism
\beq
\eU\ast A\cong \HC(A)\to\HC(A)_B\cong \eU\ast \res(A),
\eeq
where $\HC(A)_B$ is the restriction of $\HC(A)$ to $B\subset G$ extended by 0 again to $G$, and  $\res:\D_G(G)\to \D_T(T)$ is the parabolic restriction functor defined by $A\mapsto (\eU\ast A)|_T[-2\dim U](-\dim U)$ (see  (\ref{eq:resandHC})) defined by restricting the Harish-Chandra transform to the Borel $B$ and using the identification (\ref{eq:identificationBT}).

\bdefn\label{def:vanishingonG}
The category $\Dcirc_G(G)\subset \D_G(G)$ is defined to be the full subcategory formed by the objects $A\in \D_G(G)$ whose Harish-Chandra transform $e_U\ast A\cong A\ast e_U$ is supported on $B\subset G$, i.e. those $A\in \D_G(G)$ such that the natural adjunction morphism $\HC(A)=\eU\ast A\xto{}\eU\ast\res(A)$ is an isomorphism. 
\edefn

\begin{rk}\label{cent_an}
For a parabolic subgroup $P=L\cdot U_P\subset G$ one can consider 
the parabolic Harish-Chandra functor $\HC_P:A\mapsto e_{U_P}*A$.
It is not hard to see that for
$A\in \Dcirc_G(G)$ support of $HC_P(A)$ is contained in $P$. If $A$
 is equipped with the structure of a Weil sheaf, the corresponding
 trace of Frobenius function $f_A$ on the finite Chevalley group $G({\mathbb F}_q)$ also satisfies the similar vanishing condition. Let us call such a function a 
 {\em vanishing class function}.
 
 The space of vanishing class functions is a subring $Z_0$ in the center of 
 the group ring $Z\subset \Qlcl[G({\mathbb F}_q)]$. 
 It is not hard to see that an element of $Z_0$ acts by
 a scalar in every representation parabolically induced from an irreducible cuspidal representation of a Levi subgroup.
 
 We propose to view $Z_0$
 as an analogue of the Bernstein center of a $p$-adic group, the reason being
 that the similar vanishing is automatic for the distributions belonging to 
 Bernstein center. 
 \end{rk}
\begin{rk} \label{gamma_sh}
An important example of objects in $\Dcirc_G(G)$ is provided by 
the so called {\em $\gamma$-sheaves} defined by Braverman and Kazhdan
\cite{BK:03}; the fact that they lie in $\Dcirc_G(G)$ is the content of 
the {\em vanishing conjecture} stated in \cite{BK:03} and proved in 
\cite{Chen:1909} and \cite{LL:22}. 

As shown in {\em loc. cit.} the Frobenius trace of a $\gamma$-sheaf satisfies the following property: such an element in $Z$
%the center of the group ring of $G({\mathbb F}_q)$
acts by a scalar
in every Deligne-Lusztig representation. We expect that this property
holds for every Weil sheaf in $\Dcirc_G(G)$. The subalgebra in $Z$
defined by this condition can be viewed as a $G({\mathbb F}_q)$ analogue
of the {\em stable center} of a reductive $p$-adic group.
\end{rk}

\blem\label{lem:vanishingonGthick}
The full subcategory $\Dcirc_G(G)\subset \D_G(G)$ is a thick braided triangulated monoidal subcategory.
\elem
\bpf
We know that $\HC:\D_G(G)\to \eU\D_B(G)\eU$ is a triangulated monoidal functor and that the full subcategory $\eU\D_B(B)\subset \eU\D_B(G)\eU$ is a thick triangulated monoidal subcategory. Moreover, the full subcategory $\Dcirc_G(G)\subset \D_G(G)$ has been defined to be the full subcategory formed by those objects whose Harish-Chandra transform lies in this subcategory. The lemma easily follows from these observations. 
\epf

By definition, the parabolic restriction essentially agrees with the Harish-Chandra transform on $\Dcirc_G(G)$, i.e. we have $\HC(A)=\eU\ast\res(A)$ for all $A\in \Dcirc_G(G)$. Hence we obtain:
\blem\label{lem:indresmodule}
For $C\in \D_T(T),\ A\in \D_G(G)$ we have functorial morphisms $\ind(C)\ast A\to \ind(C\ast \res(A))$. Moreover, if $A\in \Dcirc_G(G)$ then the above gives us functorial isomorphisms
$$\ind(C)\ast A\xto{\cong} \ind(C\ast \res(A)).$$
\elem
\bpf
This follows immediately from (\ref{eq:avHC}), the definition of $\Dcirc_G(G)$ and that $\ind(C)=\av_B^G(e_U\ast C)$.
\epf

We have the  monoidal functors 
\beq\label{eq:HC1} \D_G(G)\xto{\HC} \Hecke\mbox{ and its restriction}\eeq
\beq\label{eq:HC2} \Dcirc_G(G)\to  e_U\D_B(B)\cong \D_T(T)\eeq
where the latter composition also agrees with the parabolic restriction functor $\res:\D_G(G)\to\D_T(T)$ on $\Dcirc_G(G)$. 

Note that $\D_B(B)\subset \D_B(G)$ is a thick braided triangulated subcategory, and moreover for any $C\in \D_B(B)$ and $D\in \D_B(G)$ we have functorial braiding isomorphisms 
\beq
\beta_{C,D}: C\ast D \xto{\cong} D*C
\eeq
which are compatible with the braiding isomorphisms in $\D_G(G)$ via the Harisch-Chandra functors (\ref{eq:HC1},\ref{eq:HC2}). In particular, the parabolic restriction functor (\ref{eq:HC2}) (or equivalently the Harish-Chandra transform) on $\Dcirc_G(G)$ is a braided monoidal functor.

For any $A_1,A_2\in \Dcirc_G(G)$ we have natural identifications
\beq\label{eq:homres}
\Hom_{\D_T(T)}(\res(A_1),\res(A_2))=\Hom_{\Hecke}(\HC(A_1),\HC(A_2))=\Hom_{\D_G(G)}(A_1,\Spr\ast A_2).
\eeq
Hence the Springer action of $W$ on $\Spr$ induces an action of $W$ on $\Hom_{\D_T(T)}(\res(A_1),\res(A_2))$.

Let $s\in W$ be a reflection corresponding to a simple root $\alpha$ with respect to a pair $T\subset B$ and let $U_{\pm\alpha}$ be the corresponding root subgroups. We have the parabolic subgroup $P_\alpha:=B\sqcup BsB$ with Levi decomposition $P_\alpha=G_\alpha U_s$, where $G_\alpha:=\<T,U_\alpha,U_{-\alpha}\>$ and $U_s:=U\cap {}^sU$. Then we have $U=U_\alpha U_s=U_sU_\alpha$ and hence  $\eU=e_{U_{\alpha}}\ast e_{U_s}=e_{U_s}\ast e_{U_\alpha}$. The Harish-Chandra transform and the adjoint averaging functors factor as 
\beq
\xymatrix@C=7em{
\D_G(G) \ar@/^/[r]^-{e_{U_s}\ast(\cdot)} & e_{U_s}\D_{P_\alpha}(G)e_{U_s}\ar@/^/[l]^{\av_{P_\alpha}^G} \ar@/^/[r]^{\eU\ast(\cdot)} & \eU\D_B(G)\eU\ar@/^/[l]^{\av_B^{P_\alpha}}.
}
\eeq
Note that for $A\in \D_G(G), C\in e_{U_s}\D_{P_\alpha}(G)e_{U_s}, D\in e_{U}\D_{B}(G)e_{U}$ we have functorial isomorphisms 
\beq
\av_{P_\alpha}^G(C\ast (e_{U_s}\ast A))\cong \av_{P_\alpha}^G(C)\ast A,
\eeq
\beq\label{eq:avHCforP}
\av^{P_\alpha}_B(D\ast (e_{U}\ast C))\cong \av^{P_\alpha}_B(D)\ast C.
\eeq
In particular, for $C\in e_{U_s}\D_{P_\alpha}(G)e_{U_s}$, $\av^{P_\alpha}_B(e_{U}\ast C)\cong \av^{P_\alpha}_B(\eU)\ast C\cong e_{U_s}\ast \Spr_{G_\alpha}\ast C$, where $\Spr_{G_\alpha}$ is the Springer sheaf on $G_\alpha$. 

For $A_1,A_2\in \D_G(G)$, adjunction gives us identifications
\beq\label{eq:adjunctions}
\begin{split}
\Hom_{\Hecke}(\eU\ast A_1,\eU\ast A_2)=\Hom_{e_{U_s}\D_{P_{\alpha}}(G)e_{U_s}}(e_{U_s}\ast A_1,\av_B^{P_\alpha}(\eU)\ast e_{U_s}\ast A_2)\\=\Hom_{\D_G(G)}(A_1,\Spr\ast A_2).
\end{split}
\eeq
\blem\label{lem:HCforP}
If $A\in \Dcirc_G(G)$, then the object $e_{U_s}\ast A\in e_{U_s}\D_{P_\alpha}(G)e_{U_s}$ is supported on $P_\alpha$ and parabolic restriction from $G$ to $G_\alpha$ defines a braided monoidal functor 
\[
\res^G_{G_\alpha}:\Dcirc_G(G)\to \Dcirc_{G_\alpha}(G_\alpha).
\]
For $A_1,A_2\in \Dcirc_G(G)$ we have natural identifications
\beq\label{eq:adjunctionsforres}\begin{split}
\Hom_{\D_T(T)}(\res(A_1),\res(A_2))=\Hom_{\D_{G_\alpha}(G_\alpha)}(\res^G_{G_\alpha}(A_1),\Spr_{G_\alpha}\ast\res^G_{G_\alpha}(A_2))\\=\Hom_{\D_G(G)}(A_1,\Spr\ast A_2).
\end{split}\eeq
\elem
\bpf
For $A\in \Dcirc_G(G)$, consider the restriction of $e_{U_s}\ast A$ to a non-trivial double coset ${P_\alpha w P_\alpha}$, namely the object $(e_{U_s}\ast A)|_{P_\alpha w P_\alpha}\in e_{U_s}\D_{P_\alpha}(G)e_{U_s}$. Since $\eU\ast A$ has support inside the Borel $B\subset P_\alpha$, we must have $\eU\ast (e_{U_s}\ast A)|_{P_\alpha w P_\alpha}=0$. Applying the functor $\av_{B}^{P_\alpha}$ and using Equation (\ref{eq:avHCforP}) we get
$$0\cong\av^{P_\alpha}_B(e_{U}\ast(e_{U_s}\ast A)|_{P_\alpha w P_\alpha})\cong \av^{P_\alpha}_B(\eU)\ast (e_{U_s}\ast A)|_{P_\alpha w P_\alpha}\cong e_{U_s}\ast \Spr_{G_\alpha}\ast (e_{U_s}\ast A)|_{P_\alpha w P_\alpha}.$$
Since $\Spr_{G_\alpha}$ contains the unit object $\delta_1$ as a direct summand, we must have $(e_{U_s}\ast A)|_{P_\alpha w P_\alpha}=0$. Hence the object $e_{U_s}\ast A$ lies in $e_{U_s}\D_{P_\alpha}(P_\alpha)\cong \D_{G_\alpha}(G_\alpha)$ and hence we have the braided monoidal functor $\res^G_{G_\alpha}:\Dcirc_G(G)\to \D_{G_\alpha}(G_\alpha)$. Moreover, $\eU\ast e_{U_s}\ast A$ has support contained inside $B$, hence we must in fact have $\res^G_{G_\alpha}(A)\in \Dcirc_{G_\alpha}(G_\alpha)$ as desired. Finally (\ref{eq:adjunctionsforres}) follows from (\ref{eq:adjunctions}).
\epf

\section{Enhanced parabolic restriction of vanishing sheaves}\label{sec:proofofthm:main1} 
We now state the following result, which contains Theorem \ref{thm:main1} as statement (iv) below: 
\bthm\label{thm:dcircg1}
(i) On the full subcategory $\Dcirc_G(G)\subset \D_G(G)$, the parabolic restriction functor (which essentially agrees with the Harish-Chandra functor on the subcategory) can be upgraded to a braided monoidal functor
\beq\label{rescirc}\res^\circ:\Dcirc_G(G)\to \D_N(T)\eeq
such that for any $A_1,A_2\in \Dcirc_G(G)$ we have natural identifications
\beq\label{eq:homrescirc}\Hom_{\D_N(T)}(\res^\circ(A_1),\res^\circ(A_2))=\Hom_{\D_T(T)}(\res(A_1),\res(A_2))^W=\Hom_{\D_G(G)}(A_1, A_2).\eeq
(ii) In fact $\res^\circ$ defines a fully faithful braided monoidal functor
\[
\res^\circ:\Dcirc_G(G)\hookrightarrow \Dcirc_N(T).
\]
(iii) Let $\Perv^\circ_G(G):=\Perv_G(G)\cap\Dcirc_G(G)$. We have the following commutative diagram of functors
\beq
\xymatrix{
\Perv^\circ_G(G)\ar[r]^-{\res^\circ}\ar@{^(->}[d] & \Perv_N(T)\ar[d]^{\cong}\\
\Perv_G(G)\ar[r]^-{\res_W} & \Perv_W(T).}
\eeq
(iv) In fact, $\res^\circ$ defines a braided monoidal triangulated equivalence
\beq
\res^\circ:\Dcirc_G(G)\xto{\cong}\Dcirc_N(T)\subset \D_N(T)
\eeq
which is also t-exact for the perverse t-structures.
\ethm

Before proving the above result, let us derive a corollary. Let $\Dcirc_N(T)\ni C\mapsto \widetilde{C}\in \Dcirc_G(G)$ be inverse to the braided equivalence $\res^\circ$. Then the composition 
\beq
\ind^\circ:\Dcirc_N(T)\subset \D_N(T) \to \D_T(T)\xto{\ind}\D_G(G)
\eeq
can be naturally identified with $C\mapsto \av_B^G(\eU\ast\widetilde{C})=\Spr\ast\widetilde{C}$. Hence we obtain
\bcor\label{cor:wactiononindcirc}
There is a natural $W$-action on the functor $\ind^\circ:\Dcirc_N(T)\to \D_G(G)$ and taking $W$-invariants we obtain a braided monoidal equivalence
\[
{\ind^\circ}^W:\Dcirc_N(T)\xto{\cong}\Dcirc_G(G)
\]
which is inverse to $\res^\circ$.
\ecor
\bpf
This follows from the Springer action of $W$ on $\Spr$ and the fact that $\Spr^W=\delta_1\in \D_G(G)$.
\epf

\subsection{Proof of Theorem \ref{thm:dcircg1}(i)}
Let us first construct the functor $\res^\circ:\Dcirc_G(G)\to \D_N(T)$. The construction of such a functor  amounts to proving that the parabolic restriction functor is independent of the choice of the Borel subgroup containing $T$ on the full subcategory $\Dcirc_G(G)\subset \DGG$.
\brk
We note that by \cite[\S1.4]{Gun:17}, parabolic restriction on the whole category $\D_G(G)$ does depend on the choice of the Borel subgroup. However we will prove that when restricted to the full subcategory $\Dcirc_G(G)$, it is independent of this choice.
\erk

Let $w\in W$ and let ${}^wB=T\cdot { }^wU$ be the Borel subgroup containing $T$ in relative position $w$ with respect to $B$. Consider the $T$-stable unipotent subgroups $U_w:= { }^wU\cap U,\  U'_w:={ }^w\Uop\cap U,\ \Uop'_w:= { }^wU\cap \Uop$. Then we have the product decompositions \[U=U_w\cdot U'_w=U'_w\cdot U_w,\ {}^wU=U_w\cdot \Uop'_w=\Uop'_w\cdot U_w\]
and natural isomorphisms of convolutions of closed idempotents
\beq
\eU\cong e_{U'_w}\ast e_{U_w}\cong e_{U_w}\ast e_{U'_w}\ \ \mbox{ and }  e_{{}^wU} \cong e_{U_w}\ast e_{\Uop'_w}\cong e_{\Uop'_w}\ast e_{U_w}.
\eeq

 Consider the triangulated category $\eU\D_T(UT\ {}^wU)e_{{ }^wU}$. We will now equip it with the structure of a braided monoidal category. Note that we have equivalences of triangulated categories 
\beq\label{eq:resB}
\begin{split}
\eU\D_T(UT\ {}^wU)e_{{ }^wU} \cong \eU\D_T(B)\cong \D_T(T)\mbox{ given by }\ \ \ \ \ \ \ \ \\ C\mapsto C|_B[-2l(w)](-l(w))\mapsto C|_T[-2\dim U-2l(w)](-\dim U -l(w)),
\end{split}
\eeq
\beq\label{eq:reswB}
\begin{split}
\eU\D_T(UT\ {}^wU)e_{{ }^wU} \cong \D_T({}^wB)e_{{ }^wU}\cong \D_T(T)\mbox{ given by }\ \ \ \ \ \ \ \ \\ C\mapsto C|_{{}^wB}[-2l(w)](-l(w))\mapsto C|_T[-2\dim U-2l(w)](-\dim U -l(w)).
\end{split}
\eeq

Now each of $\D_T({}^wB)e_{{ }^wU}=\D_{{}^wB}({}^wB)e_{{ }^wU}$ and $\eU\D_T(B)=\eU\D_B(B)$ are braided monoidal categories which are both braided monoidally equivalent to $\D_T(T)$. Hence using the above triangulated equivalence, we can equip $\eU\D_T(UT\ {}^wU)e_{{ }^wU}$  with the structure of a braided monoidal category using that of $\D_T(T)$. 

Let $B'$ be any Borel subgroup of $G$ containing $T$ and let $U'$ be its unipotent radical. Consider the parabolic restriction functor with respect to $B'$ (see also Equation \ref{eq:resandHC}),
$\res_{B'}:\D_G(G)\to \D_T(T)\mbox{ defined by }\D_G(G)\ni A\mapsto (e_{U'}\ast A)|_{T}[-2\dim U](-\dim U).$ 
If $A\in \Dcirc_G(G)$, the support of $e_{U'}\ast A$ is contained inside $B'$ and we have the braided monoidal functor 
\beq
\res_{B'}:\Dcirc_G(G)\to \eU\D_{B'}(B')\cong \D_T(T) \mbox{ defined by }A\mapsto (e_{U'}\ast A)|_{T}[-2\dim U](-\dim U).
\eeq
Now we consider the two Borel subgroups $B$ and ${ }^wB$ which contain $T$ and consider the composition
\beq\begin{split}
\Dcirc_G(G)\rar{}\eU\D_T(UT{ }^wU)e_{{ }^wU}\cong \D_T(T),\mbox{ defined by }\ \ \ \ \ \ \ \ \\
 A\mapsto (\eU\ast A\ast e_{{}^wU})|_T[-2\dim U-2l(w))](-\dim U-l(w)).
 \end{split}
\eeq
Computing the second functor in the above composition using the Borel $B$ as in (\ref{eq:resB}) we see that the composition above is naturally isomorphic to $\res=\res_B:\Dcirc_G(G)\to \D_T(T)$. On the other hand, using the Borel ${ }^wB$ as in (\ref{eq:reswB}) we see that it is also naturally isomorphic to $\res_{{}^wB} = w\circ \res:\Dcirc_G(G)\to \D_T(T)$.
This gives us a natural isomorphism between the two braided monoidal triangulated functors \beq\label{eq:wequiv}\res\cong w\circ\res:\Dcirc_G(G)\to \D_T(T)\eeq as desired. Hence we obtain a braided monoidal triangulated functor
\beq
\res^\circ:\Dcirc_G(G)\to \D_N(T)=(\D_T(T))^W.
\eeq
This induces an action of $W$ on the morphism spaces $\Hom_{\D_T(T)}(\res(C),\res(D))$ for each $C,D\in \Dcirc_G(G)$. Namely for each $w\in W$ use (\ref{eq:waction}) to define the action of $w$:
\beq\label{eq:waction} w:\left(\res(C)\rar{f}\res(D)\right)\mapsto \left(\res(C)\cong w\circ\res(C)\xto{w(f)}w\circ\res(D)\cong \res(D)\right).\eeq
For $C,D\in \Dcirc_G(G)$ we also have the identification  (\ref{eq:homres}) obtained using adjunction
\beq\label{eq:springeradjunction}\Hom_{\D_T(T)}(\res(C),\res(D))=\Hom_{\D_G(G)}(C,\Spr\ast D).\eeq
We have a $W$-action on the right hand side of (\ref{eq:springeradjunction}) coming from the $W$-action on the Springer sheaf $\Spr$. We will now prove that these $W$-actions match. It is enough to consider the action of the simple reflections. Let $s\in W$ be a reflection corresponding to a simple root $\alpha$, $U_{\pm\alpha}$  the corresponding root subgroups and $P_\alpha=B\sqcup BsB$ the parabolic subgroup with Levi decomposition $P_\alpha=G_\alpha U_s$, where $G_\alpha:=\<T,U_\alpha,U_{-\alpha}\>$ and $U_s:=U\cap {}^sU$. Using Lemma \ref{lem:HCforP}, to compare the action of $s$ coming from the above construction with that coming from the Springer action it is enough to restrict our attention to the subgroup $G_\alpha$ of semisimple rank 1. Hence let us assume without loss of generality that $G$ is of semisimple rank 1. In this case $U,\Uop\cong \Ga$, with $T$ acting via the characters $\alpha, \alpha^{-1}$ respectively and $W=\<s\>$ is cyclic of order 2.

By functoriality it is enough to prove that the actions of $s$ agree when $D=\delta_1$, the unit object of $\Dcirc_G(G)$, namely we will compare the actions of $s$ on:
\beq\label{eq:tocomparewactions}
\Hom_{\D_T(T)}(\res(C),\delta_1)=\Hom_{\D_G(G)}(C,\Spr).
\eeq

Also in this case the Springer sheaf splits as $\Spr=\delta_1\oplus \Pcal_\Ucal$, where $\Pcal_{\Ucal}:={\Qlcl}_{\Ucal}[2](1)$ is the constant perverse sheaf supported on the unipotent cone $\Ucal\subset G$. The adjunction morphism $\Spr\to \eU$ in $\D_B(G)\subset \D_T(G)$ corresponds to the natural morphisms $\delta_1\to \eU$ and $\Pcal_{\Ucal}\to \eU$, the later defined from the closed immersion $U\subset \Ucal$. The group $\<s\>$ acts trivially on $\delta_1$ and by the sign on $\Pcal_{\Ucal}$. In this case Equation (\ref{eq:tocomparewactions}) becomes
\beq\label{eq:tocomparesactions}
\Hom_{\D_T(T)}(\res(C),\delta_1)=\Hom_{\D_G(G)}(C,\delta_1)\oplus \Hom_{\D_G(G)}(C,\Pcal_{\Ucal}).
\eeq

Recall that $\eU, \eUop\in \D_T(G)$ are closed idempotents.
Hence by \cite[\S2.8]{BoDr:14}, for each $C\in \Dcirc_G(G)$ we have natural identifications
\[\Hom_{\D_T(T)}(\res(C),\delta_1)=\Hom_{\D_T(G)}(\eU\ast C\ast\eUop,\eU\ast\eUop)=\Hom_{\D_T(G)}(\eU\ast C,\eU\ast \eUop).\]
It will be more convenient for us to consider the action of $s$ on the space $\Hom_{\D_T(G)}(\eU\ast C,\eU\ast \eUop)\cong \Hom_{\D_G(G)}(C,\Spr)$. To prove the compatibility of the $s$-actions, we must prove that $s$ acts trivially on the image of $\Hom_{\D_G(G)}(C,\delta_1)\hookrightarrow \Hom_{\D_T(G)}(\eU\ast C,\eU\ast \eUop)$ and by $-1$ on the image of $\Hom_{\D_G(G)}(C,\Pcal_{\Ucal})\hookrightarrow \Hom_{\D_T(G)}(\eU\ast C,\eU\ast \eUop)$. 

The former assertion follows from the commutative diagram
\[
\xymatrix{
\delta_1\ar[r]\ar[d] & \eU\ar[d]\\
\eUop\ar[r] & \eU\ast \eUop
}
\]
in $\D_T(G)$. For the latter: Consider a morphism $C\xto{f} {\Pcal}_{\Ucal}$ in $\D_G(G)$. Note that we have the two natural morphisms $\Pcal_{\Ucal}\xto{\alpha} \eU$ and $\Pcal_{\Ucal}\xto{\overline{\alpha}} \eUop$. Then we need to compare the following two morphisms in $\D_T(G)$ and prove that they are negatives of each other:
\beq\label{eq:twomorphisms}
\begin{split}
    \eU\ast C\xto{\eU\ast f} \eU\ast\Pcal_{\Ucal}\xto{\eU\ast\alpha} \eU\to \eU\ast\eUop\\
    \eU\ast C\xto{\eU\ast f} \eU\ast\Pcal_{\Ucal}\xto{\eU\ast\overline{\alpha}} \eU\ast\eUop.
\end{split}
\eeq
Since $\eU\ast C$ is supported on $B$, it is enough to consider the restriction of the above morphisms to $B$. Consider the following two compositions in $\D_T(G)$ which come from the two $T$-equivariant inclusions $\{1\}\subset U$ and $\{1\}\subset \Uop$ respectively:
\beq
\begin{split}
    \delta_1\to \eU \to \delta_1[2](1)\\
    \delta_1\to \eUop \to \delta_1[2](1).
\end{split}
\eeq
Since $T$ acts on $U$ and $\Uop$ by opposite characters, the above compositions correspond to negatives of each other in the $T$-equivariant cohomology $\RHom_{\D_T(1)}(\Qlcl,\Qlcl)=\RHom_{\D_T(G)}(\delta_1,\delta_1)$. From this it follows that the restriction to $B$ of the two morphisms in Equation \ref{eq:twomorphisms} are negatives of each other as desired.

Hence we have proved that the identification in Equation \ref{eq:springeradjunction} is $W$-equivariant. Taking $W$-invariants and noting that $\Spr^W=\delta_1$, we obtain
\beq
\begin{split}
\Hom_{\D_N(T)}(\res^\circ(C),\res^\circ(D))=\Hom_{\D_T(T)}(\res(C),\res(D))^W\\
=\Hom_{\D_G(G)}(C,\Spr\ast D)^W=\Hom_{\D_G(G)}(C,D)
\end{split}
\eeq
completing the proof of Theorem \ref{thm:dcircg1}(i).

\subsection{Proof of Theorem \ref{thm:dcircg1}(ii)}\label{sec:proofothm:dcircg1:2}
The fact that the functor $\res^\circ:\Dcirc_G(G)\to \D_N(T)$ is fully faithful is just (\ref{eq:homrescirc}) in Theorem \ref{thm:dcircg1}(i). Let us now prove that the essential image is contained inside $\Dcirc_N(T)$. 

We must show that for each $A\in \Dcirc_G(G)$ the object $\res^\circ(A)\in \D_N(T)$ is central on $T$, i.e. lies in $\Dcirc_N(T)$. Let $s\in W$ be a simple reflection corresponding to a simple root $\alpha \in \Delta$ and corresponding Levi subgroup $G_{\alpha}$. Since the parabolic restriction of $A$ to $G_{\alpha}$ lies in $\Dcirc_{G_{\alpha}}(G_{\alpha})$, by Lemma \ref{lem:centralconditions}(iv) we may assume without loss of generality that $G=G_{\alpha}$, i.e. $G$ is of semisimple rank 1. 

In this case, it is straightforward to check that a multiplicative local system $\chi\in \Ccal(T)$ can be extended to the whole of $G$ if and only if ${\alpha^\vee}^*\chi\cong \Qlcl$. For such a $\chi$, let the multiplicative local system $\chi_G$ denote its extension to $G$. Then $\chi_G\in \D_G(G)$. Note that $\chi_G|_{U}\cong \Qlcl$ and hence the restriction of $\chi_G$ to each coset $Ug$ is also constant. It then follows that $A\otimes\chi_G\in \Dcirc_G(G)$. It also follows from the connectedness of $G$ that the adjoint action of every element $g\in G$ on the global cohomology $H^*_c(G,A\otimes \chi_G)$ is trivial. On the other hand we have 
\[H^*_c(G,A\otimes \chi_G)=H^*_c(G,\eU\ast(A\otimes\chi_G))=H^*_c(B,\eU\ast(A\otimes \chi_G))=H^*_c(T,\res(A)\otimes \chi).\]
Hence the induced action of $s\in W$ on $H^*_c(T,\res^\circ(A)\otimes \chi)$ is trivial for each $\chi\in\Ccal(T)$ such that ${\alpha^\vee}^*\chi\cong \Qlcl$. Hence $\res^\circ(A)\in \Dcirc_{N}(T)$ as desired.

\subsection{Proof of Theorem \ref{thm:dcircg1}(iii)}\label{sec:dgcirc1iii}
To prove Theorem \ref{thm:dcircg1}(iii), we will use the following result, namely the vanishing conjecture proved in \cite[Thm. 1.5]{Chen:1909}, see also \cite{BITV}:
\bthm\label{thm:van}
The $W$-invariant induction functor $\ind^W$ takes a $W$-equivariant central complex on $T$ to a vanishing complex in $\Dcirc_G(G)$, i.e. we have the functor
\[
\ind^W:\Dcirc_W(T)\to \Dcirc_G(G).
\]
\ethm
In \cite{Chen:1909} this result is proved for $p$ not dividing some integer $M$ depending on the type of $G$. However this restriction can be removed as follows. For unipotently monodromic perverse sheaves in $\Dcirc_W(T)$, this is proved in \cite[Cor. 5.7.2]{BITV}. In a work under preparation by the same authors, this result will be extended to monodromic perverse sheaves in $\Dcirc_W(T)$ for any monodromy (corresponding to any $W$-orbit in $\Ccal(T)$). Finally, as proved in \cite{Chen:1909}, the vanishing conjecture follows for any object in $\Dcirc_W(T)$ if it is known for all monodromic perverse objects. We also remark that in the proof of Theorem \ref{thm:dcircg1}(iii) below, we will only use the above vanishing result for monodromic perverse sheaves in $\Dcirc_W(T)$.

Since parabolic induction and restriction are t-exact for the perverse t-structure by \cite{BYD:18}, we have a fully faithful functor
\beq
\res^\circ: \Perv^\circ_G(G)\to \Perv^\circ_N(T)\cong\Perv^\circ_W(T)
\eeq
and the functor
\beq
\ind^W:\Perv^\circ_W(T)\to \Perv^\circ_G(G).
\eeq
Note that on the category $\Perv^\circ_G(G)$ we have the two functors (from Definition \ref{defn:resw} and Theorem \ref{thm:dcircg1}(i)) 
\[\res_W,  \res^\circ:\Perv^\circ_G(G)\to \Perv_W(T)\]
both of whose compositions with the forgetful functor $\Perv_W(T)\xto{\forg_W} \Perv(T)$ give the parabolic restriction functor: 
\beq
\res=\forg_W\circ\res^\circ=\forg_W\circ\res_W: \Perv^\circ_G(G)\to \Perv(T).
\eeq
Also, by Lemma \ref{lem:resWindW}(ii) we know that 
\beq
\res_W\circ\ind^W=\Id_{\Perv^\circ_W(T)}.
\eeq 
Combining the above observations we have
\beq\label{eq:forgetweq}
\forg_W\circ \res^\circ\circ\ind^W=\forg_W:\Perv^\circ_W(T)\to \Perv(T).
\eeq
We will now prove that the composition
\beq\label{eq:reswindw=id}
\res^\circ\circ\ind^W: \Perv^\circ_W(T)\to \Perv^\circ_W(T)
\eeq
can be naturally identified with the identity functor on $\Perv_W^\circ(T)$. 
% It follows from (\ref{eq:forgetweq}) that 
% \beq\label{eq:forgetwderived}
% \forg_W\circ\res^\circ\circ\ind^W=\forg_W: \D^\circ_W(T)\to \D(T).
% \eeq 
By \ref{eq:forgetweq} it remains to check that two specific $W$-equivariance structures defined on the same underlying object of $\Perv(T)$ agree.

Let us first consider objects $M\in \Perv^\circ_W(T)\subset \Dcirc_W(T)$ which are shifts of tame $W$-equivariant local systems on $T$, namely the monodromic perverse sheaves in $\Dcirc_W(T)$.  Such local systems correspond to $\Qlcl$-representations (say $V$) of $W\ltimes \pi_1^t(T)$. On the other hand, we have a surjective homomorphism $\pi_1(G^{rs})\to W\ltimes \pi_1^t(T)$ from the \'{e}tale fundamental group of the regular semisimple locus $G^{rs}\subset G$. For such an $M\in \Perv^\circ_W(T)$ corresponding to the representation $V$, we have $\ind^WM=\IC(G^{rs},V)\in \Perv^\circ_G(G)$. Let $T^\reg\subset T$ denote the open subset of regular elements.

The orbits for the conjugation action of $U$ on $T^\reg U$ are just the cosets $tU=Ut$ for $t\in T^\reg$. Hence we see that $\eU\ast (\ind^W M)|_{T^\reg U}=(\ind^WM)|_{T^\reg U}$. It follows that (see also Equation \ref{eq:resandHC})
\beq
(\res\circ\ind^WM)|_{T^\reg}=(\ind^WM)|_{T^\reg}[-2\dim U](-\dim U)=M|_{T^\reg}.
\eeq
The same is true for the parabolic restriction $\res_{U'}$ with respect to any maximal unipotent $U'$ normalized by $T$. Hence it follows that $(\res^\circ\circ\ind^WM)|_{T^\reg}=M|_{T^\reg}$. By (\ref{eq:forgetweq}), we know that the underlying object in $\Perv(T)$ corresponding to $\res^\circ\circ\ind^WM\in \Perv^\circ_W(T)$ is the same as the underlying perverse sheaf corresponding to $M$ itself, which is just a shifted local system. Moreover, we have seen above that the $W$-equivariance structures on the shifted local system $\forg_WM$ on $T$ corresponding to $\res^\circ\circ\ind^WM$ and $M$ agree when restricted to the dense open $T^\reg\subset T$. Hence we see that $\res^\circ\circ\ind^WM=M$ for all $M\in \Perv^\circ_W(T)$ of the type above. Hence the statement also follows for objects of the full triangulated subcategory of $\Dcirc_W(T)$ generated by such objects, namely all the monodromic $W$-equivariant central complexes on $T$.

Next we consider an arbitrary object $(C,\varphi_C)\in \Perv^\circ_W(T)$, where we specify the $W$-equivariance structure on $C\in \D(T)$ by $\varphi_C$. Namely for every $w\in W$, we have the isomorphism $\varphi_C(w):w^*C\xto{\cong} C$.  By Equation \ref{eq:forgetweq} we know that $\res^\circ\circ\ind^W(C,\varphi_C)=(C,\psi_C)$, where $\psi_C$ is some $W$-equivariance structure on $C\in \D(T)$. It remains to show that $\varphi_C=\psi_C$.

Consider a morphism $(C,\varphi_C)\xto{f}(M,\varphi_M)$ in $\D^\circ_W(T)$, where $(M,\varphi_M)\in \D^\circ_W(T)$ is monodromic as in the previous paragraph.  For these objects we have proved that $\res^\circ\circ\ind^W(M,\varphi_M)=(M,\varphi_M)$.
Hence applying $\res^\circ\circ\ind^W $ to the morphism $f$, we obtain the commutative diagram
\beq\label{eq:diagramCtoM}
\xymatrix@C=5.5em{
  C \ar[r]^{\psi_C(w)^{-1}} \ar[d]_{f} & w^*C \ar[r]^{\varphi_C(w)} \ar[d]^{w^*f} & C \ar[d]^f \\
  M \ar[r]^{\varphi_M(w)^{-1}}\ar@{ = }@/_1pc/[rr] & w^*M \ar[r]^{\varphi_M(w)} & M \\
}
\eeq
\noindent for all morphisms $f$ as above. We will now use this diagram to deduce that for each $w\in W$, the composition $\varphi_C(w)\circ\psi_C(w)^{-1}$ along the top row of the diagram equals the identity of $C$. Let us denote $\xi_{C,w}:=\varphi_C(w)\circ\psi_C(w)^{-1}$.

We will use the Mellin transform $\Mcal_!:\D(T)\to D^b\Coh\Ccal(T)$, which is a  t-exact functor (see \cite[Thm. 3.4.1]{GL}), where $\D(T)$ is equipped with the perverse t-structure and $D^b\Coh\Ccal(T)$ is equipped with the t-structure dual to its standard t-structure. Also by \cite[Prop. 3.4.6]{GL}, $\Mcal_!A=0$ if and only if $A=0$. From these two observations it follows that it is enough to prove that the Mellin transform of $\xi_{C,w}$ is the identity on $\Mcal_!C$.

We also recall that the above Mellin transform gives an equivalence between the full subcategories of monodromic complexes in $\D(T)$ and coherent complexes with finite set-theoretic support in $D^b\Coh(\Ccal(T))$. Moreover, for $C,M\in \D(T)$, $\Hom(C,M)$ maps isomorphically to $\Hom(\Mcal_!C,\Mcal_!M)$ when $M$ is monodromic, this follows from \cite[Cor. 3.3.2]{GL}.

Let $\Theta\subset \Ccal(T)$ be the $W$-orbit of a multiplicative local system $\chi\in \Ccal(T)$. We have the quotient map $\Ccal(T)\xto{\pi} \Ccal(T)//W_\chi$. For a positive integer $n$, let $\mathbb{O}^{(n)}_\chi\in D^b\Coh_\chi\Ccal(T)$ be the structure sheaf of the inverse image $\mathbb{D}^{(n)}_\chi\subset \Ccal(T)$ under $\pi$ of the $n$-th nilpotent thickening of the closed point $\pi(\chi)$ in $\Ccal(T)//W_\chi$. Then $\mathbb{D}^{(n)}_\chi$ is some nilpotent thickening of $\chi\in \Ccal(T)$ and let $\mathbb{D}_\Theta^{(n)}:=W\cdot\mathbb{D}^{(n)}_\chi\subset \Ccal(T)$, the corresponding nilpotent thickening of the orbit $\Theta$. Let $\mathbb{O}^{(n)}_\Theta\in D^b\Coh_\Theta\Ccal(T)$ its structure sheaf, which is set theoretically supported on the finite set $\Theta\subset \Ccal(T)$. It has a natural $W$-equivariance structure. Also since $\mathbb{O}^{(n)}_\chi$ comes via pullback from the quotient map $\pi$, the inverse Mellin transform $\bbe^{(n)}_{\Theta}$ of $\mathbb{O}^{(n)}_\Theta$ lies in $\D^\circ_W(T)$. We have the natural morphism $\delta_1\to \bbe^{(n)}_{\Theta}$ in $\D^\circ_W(T)$, whose Mellin transform is the natural projection $\Ocal_{\Ccal(T)}\onto \mathbb{O}^{(n)}_{\Theta}$ from the structure sheaf of $\Ccal(T)$ to the structure sheaf 
 the closed subscheme $\iota:\mathbb{D}^{(n)}_{\Theta}\hookrightarrow \Ccal(T)$. Now convolving this morphism with $C$ we get a morphism $C\xto{f} \bbe^{(n)}_{\Theta}\ast C$ in $\D^\circ_W(T)$ whose Mellin transform is the natural restriction morphism $\Mcal_!C\to \iota_*\iota^*(\Mcal_!C)$. Note that the objects $\bbe^{(n)}_{\Theta}$ and the convolution $\bbe^{(n)}_{\Theta}\ast C\in \D^\circ_W(T)$ are monodromic. 
 
 Hence applying the Mellin transform to the diagram (\ref{eq:diagramCtoM}) for the arrow $C\xto{f} \bbe^{(n)}_{\Theta}\ast C$ we obtain the commutative diagram:
 \beq\label{eq:Mellindiagram}
\xymatrix@C=5.5em{
\Mcal_!C \ar[r]^{\Mcal_!\xi_{C,w}} \ar[d]  & \Mcal_!C \ar[d] \\
  \iota_*\iota^*(\Mcal_!C) \ar@{=}[r] & \iota_*\iota^*(\Mcal_!C). \\
}
\eeq
By the adjunction $\Hom(\Mcal_!C,\iota_*\iota^*(\Mcal_!C))=\Hom(\iota^*(\Mcal_!C),\iota^*(\Mcal_!C))$, the commutativity of the above diagram is equivalent to the equality $\iota^*(\Mcal_!\xi_{C,w})=\id_{\iota^*(\Mcal_!C)}$. Hence this equality holds for all the inclusions $\iota:\mathbb{D}^{(n)}_\Theta\hookrightarrow \Ccal(T)$ of the nilpotent thickenings of all $W$-orbits in $\Ccal(T)$, which is a (typically infinite) disjoint union of Noetherian regular affine schemes. We will now use Krull's intersection theorem to conclude that $\Mcal_!\xi_{C,w}-\id_{\Mcal_!C}=0$. To see this, let $R$ be a commutative Noetherian ring and $a:M_1\to M_2$ be a morphism in $D^b\Coh(R)$. Then we have $\Hom(M_1,M_2)=\Hom(R,\mathcal{RH}{om}(M_1,M_2))=H^0(\mathcal{RH}om(M_1,M_2))$. Then using Krull's intersection theorem for the finitely generated $R$-module $H^0(\mathcal{RH}om(M_1,M_2))$, it follows that the morphism $a$ is 0 if and only if its restriction to each nilpotent thickening in $\operatorname{Spec}(R)$ is zero.

As we noted before, $\Mcal_!\xi_{C,w}=\id_{\Mcal_!C}$ implies that $\xi_{C,w}=\id_C$ as desired. % and we have already proved the result for such objects. Hence $\bbe^{(n)}_{\Theta}\ast \varphi_C(w)=\bbe^{(n)}_{\Theta}\ast \psi_C(w)$, i.e. $\bbe^{(n)}_{\Theta}\ast f_{C,w}$ is the identity on $\bbe^{(n)}_{\Theta}\ast C$. Applying the Mellin transform, the convolution goes to the (derived) tensor product and we have $\mathbb{O}^{(n)}_\Theta\otimes \Mcal_!f_{C,w}$ is the identity on $\mathbb{O}^{(n)}_\Theta\otimes \Mcal_!C$ for every $W$-orbit $\Theta\subset \Ccal(T)$ and for every positive integer $n$. In other words, the (derived) restriction of the morphism $\Mcal_!f_{C,w}:\Mcal_!C\to \Mcal_!C$ to each nilpotent thickening $\mathbb{D}^{(n)}_\Theta\subset\Ccal(T)$ of each $W$-orbit $\Theta\subset \Ccal(T)$ is the identity. Hence $\Mcal_!f_{C,w}$ is the identity, and thus $f_{C,w}$ is the identity as desired.
Hence  we have shown that $\res^\circ\circ\ind^W$ can be identified with the identity functor on $\Perv^\circ_W(T)$. On the other hand by Theorem \ref{thm:dcircg1}(i),(ii) we have seen that $\res^\circ:\Perv^\circ_G(G)\to \Perv^\circ_W(T)\subset \Perv_W(T)$ is fully faithful. Hence $\res^\circ$ and $\ind^W$ define inverse equivalences between the categories $\Perv^\circ_G(G)$ and $\Perv^\circ_W(T)$. Comparing with Lemma \ref{lem:resWindW}(ii) gives us an identification between the two functors $\res^\circ=\res_W:\Perv^\circ_G(G)\xto{\cong} \Perv^\circ_W(T)$ as desired.

% Let $A\in \Perv^\circ_G(G)\subset \Dcirc_G(G)$. Then $\res(A)\in \Perv(T)$ has two $W$-equivariant structures coming from $\res_W(A)$ and $\res^\circ(A)$. In other words for each $w\in W$, we have two isomorphisms between the objects $\res(A)$ and $w^*\res(A)$ in $\Perv(T)$ and we will show that they are equal.

\subsection{Proof of Theorem \ref{thm:dcircg1}(iv)}\label{sec:dgcirc1iv}
We have already seen that $\res^\circ:\Dcirc_G(G)\to \Dcirc_N(T)$ is a triangulated fully faithful braided monoidal functor. To complete the proof we must prove essential surjectivity. Let $C\in \Perv^\circ_N(T)=\Perv^\circ_W(T)$. Then by the previous subsection we have that $\res^\circ\circ\ind^W(C)=C$. Hence $\Perv^\circ_N(T)\subset \Dcirc_N(T)$ lies in the essential image. But since the thick triangulated subcategory generated by $\Perv^\circ_N(T)$ is the whole of $\Dcirc_N(T)$, we see that $\res^\circ:\Dcirc_G(G)\to \Dcirc_N(T)$ is an equivalence.

\section{The category $\D(T)$ as a bimodule category}\label{sec:bimodulecategory}
In this section we define and study the structure of a $\eUop\D(G)\eUop-\Wh$-bimodule category on $\D(T)$. We will then use this structure to study the bi-Whittaker category $\Wh$ and  complete the proof of Theorem \ref{thm:main2}. Let begin with a study of the Yokonuma-Hecke category $\Yok$.

\brk\label{rk:integrals}
In this section as well as the next one, we will often use the following notations and conventions in computations involving tensor products, $!$-pushforwards and $*$-pullbacks of $\Qlcl$-complexes. For an $\Fcal\in \D(X)$ and $x\in X$, by the `fiber' or `stalk' of $\Fcal$ at $x$ we mean the $*$-pullback of $\Fcal$ to $x\in X$, and denote it as $\Fcal(x)$. Let $f:X\to Y$ be a map of schemes and let $\Fcal\in \D(X)$. Then the notation `$\Gcal(y)=\int\limits_{\substack{x\in X\\y=f(x)}}\Fcal(x) $ for $y\in Y$' simply stands for $\Gcal=f_!\Fcal$. We note that the proper base change theorem for $\Qlcl$-complexes is implicit in this notation. If we have an isomorphism $X'\xto{\phi} X$ of schemes and an isomorphism $\Fcal'\cong \phi^*\Fcal$ in $\D(X')$, then the corresponding isomorphism $f_!\Fcal\cong f'_!\Fcal'$ (where $f':= f\circ \phi:X'\to Y$) will be denoted as 
$$\int\limits_{\substack{x\in X\\y=f(x)}}\Fcal(x)\cong \int\limits_{\substack{x'\in X'\\y=f'(x')}}\Fcal'(x').$$  With this notation, the projection formula becomes $\int\limits_{\substack{x\in X\\y=f(x)}}\Gcal(y)\otimes\Fcal(x)\cong \Gcal(y)\otimes \int\limits_{\substack{x\in X\\y=f(x)}}\Fcal(x)$. It is important to note that although the isomorphisms involving the `$\int$' above have been written at the level of stalks, the isomorphisms themselves are globally defined. We use this notation for convenience in long computations to avoid introducing new notations for schemes and morphisms that come up in the proofs.
\erk

\subsection{The Yokonuma-Hecke category}\label{sec:yokonuma}
Let us recall the construction of some special (shifted) perverse sheaves in the Yokonuma-Hecke category $\eUop\D(G)\eUop$ parametrized by elements of the Weyl group $W$. These perverse sheaves were studied in \cite{KazLau} and are the geometric analogues of the generators of the Yokonuma-Hecke algebra constructed by Juyumaya in \cite{Ju:98}. 

\subsubsection{The case of $SL_2$}
We first carry out the construction for $G=SL_2$. Note that the self convolution of the constant sheaf on $SL_2$ is given by \beq\label{eq:eSL2*eSL2}
\Qlcl\ast\Qlcl=H^*_c(SL_2,\Qlcl)\otimes \Qlcl
\eeq and the object 
\beq\label{eq:esl2}
e_{SL_2}:=\Qlcl[6](3)
\eeq is a quasi-idempotent in $\eUop\D(SL_2)\eUop$.

For each $x\in \Ga\cong \Ga^*$ we have the corresponding multiplicative local system $\L_x$ on $\Ga$. It is easy to check that we have a $\Uop=U^-$-biequivariant morphism $\pi_{12}:SL_2\to \Ga$ defined by $\begin{pmatrix}a & b\\  c & d\end{pmatrix}\mapsto b$. Define 
$$\eUop\D(SL_2)\eUop\ni \Kcal_x:=\pi_{12}^*\L_x[4](2) \mbox{ and we set }\Kcal:=\Kcal_{1}=\pi_{12}^*\L_\psi[4](2).$$
Note that for $t\in \Gm\cong T\subset SL_2$ we have ${}^t\Kcal_x=\ad(t)^*\Kcal_x\cong \Kcal_{t^2x}$.

Next, let us compute the convolution $\Kcal\ast\Kcal$ on $SL_2$. We see that (see Remark \ref{rk:integrals})
\beq\label{eq:KK1}\Kcal\ast \Kcal\begin{pmatrix}a & b\\  c & d\end{pmatrix}=\int\limits_{\tiny{\begin{pmatrix}x\hspace{-7pt} & y\\ z\hspace{-7pt} & w\end{pmatrix}}\in SL_2}\L_\psi(bx-ay)\otimes\L_\psi(y)[8](4)=\int\limits_{\tiny{(x,y)\in \bbA^2\setminus \{0\}}}\L_\psi(bx+(1-a)y)[6](3).\eeq
We have $\Uop=\left\{\begin{pmatrix}a & b\\  c & d\end{pmatrix}\in SL_2|b=1-a=0\right\}$. Hence using the cohomology vanishing of non-trivial multiplicative local systems on $\bbA^2\cong \Ga^2$
\beq\label{eq:KK2}\int\limits_{\tiny{(x,y)\in \bbA^2}}\L_\psi(bx+(1-a)y)[6](3)=e_{\Uop}\begin{pmatrix}a & b\\  c & d\end{pmatrix}.\eeq
We also have
\beq\label{eq:KK3}\int\limits_{\tiny{(x,y)=0\in \bbA^2}}\L_\psi(bx+(1-a)y)[6](3)=\Qlcl[6](3).\eeq
Now using Equations \ref{eq:KK1}, \ref{eq:KK2}, \ref{eq:KK3} and the triple $\bbA^2\setminus \{0\}\subset \bbA^2\supset \{0\}$ we obtain a distinguished triangle
\beq\label{eq:J*J}
\Kcal\ast\Kcal\to \eUop \to \Qlcl[6](3)=e_{SL_2}\to.
\eeq
\subsubsection{General reductive groups}\label{sec:sheavesinyokonuma}
Let $\alpha:T\to \Gm$ be a simple root of $G$. Let us fix a pinning, i.e. a non-degenerate $\L\in U^*$. Hence we have the homomorphism $f_\alpha:SL_2\to G$. We define 
$$\Kcal_{s_\alpha}:=\eUop\ast{f_\alpha}_!\Kcal\ast\eUop\in \eUop\D(G)\eUop,$$
$$e_{s_\alpha}:=\eUop\ast{f_\alpha}_!e_{SL_2}\ast\eUop\in \eUop\D(G)\eUop.$$
Observe that the support  of both $\Kcal_{s_\alpha}$ and $e_{s_\alpha}$ is contained inside the parabolic subgroup $\Pop_\alpha=\Bop\sqcup\Bop s_\alpha\Bop$. Let $G_\alpha$ be the corresponding Levi subgroup. We have the Levi decomposition $\Pop_\alpha=G_\alpha\Uop_{s_{\alpha}}$ where for any $w\in W$, $\Uop_w:=\Uop\cap {}^w\Uop$. The derived subgroup $G_\alpha'\subset G_\alpha=\langle T,U_\alpha,U_{-\alpha}\rangle$ equals the image $f_\alpha(SL_2)$. We see that the support of both $\Kcal_{s_\alpha},e_{s_\alpha}$ equals $\Pop'_\alpha:=G_\alpha'\Uop_{s_{\alpha}}$ where both are just shifted local systems.
We have
\blem\label{lem:equivofmonoidalcat}
The pair $(\Uop_{s_\alpha},\Qlcl)$ is a Heisenberg admissible pair for $\Pop_\alpha$. We have the equivalences of triangulated monoidal  categories $$e_{\Uop_{s_\alpha}}\D(\Pop_\alpha)\cong \D(G_\alpha)\mbox{ and }e_{\Uop_{s_\alpha}}\D(\Pop'_\alpha)\cong \D(G'_\alpha),\mbox{ as well as}$$ $$\eUop\D(\Pop_\alpha)\eUop\cong e_{U_{-\alpha}}\D(G_\alpha)e_{U_{-\alpha}}\mbox{ and }\eUop\D(\Pop'_\alpha)\eUop\cong e_{U_{-\alpha}}\D(G'_\alpha)e_{U_{-\alpha}}.$$ We have 
$$\Kcal_{s_\alpha}\cong e_{\Uop_{s_\alpha}}\ast {f_{\alpha}}_!\Kcal\cong {f_{\alpha}}_!\Kcal\ast e_{\Uop_{s_\alpha}}\mbox{ and }$$
$$e_{s_\alpha}\cong e_{\Uop_{s_\alpha}}\ast {f_{\alpha}}_!e_{SL_2}\cong {f_{\alpha}}_!e_{SL_2}\ast e_{\Uop_{s_\alpha}}.$$
Moreover $\Kcal_{s_\alpha}[-\dim \Uop]$ and $e_{s_\alpha}[-\dim \Pop'_\alpha]$ are local systems supported on $\Pop_\alpha'$. The ranks of both these local systems equal 1 in case $G_\alpha'\cong SL_2$ and 2 in case $G'_\alpha\cong PGL_2$.
\elem

\blem\label{lem:qiesalpha}
(i) For each simple root $\alpha$, the object $e_{s_\alpha}$ is a quasi-idempotent in $\eUop\D(G)\eUop$, and we have 
$$e_{s_\alpha}\ast e_{s_\alpha}\cong H^*_c(SL_2,\Qlcl)[6](3)\otimes e_{s_\alpha}.$$\\
(ii) For each simple root $\alpha$, we have the distinguished triangle
$$\Kcal_{s_\alpha}\ast \Kcal_{s_\alpha}\to \eUop\to e_{s_\alpha}\to$$
\elem
\bpf
Both the statements follow easily using the fact that the functor ${f_\alpha}_!$ is monoidal, using Equations (\ref{eq:eSL2*eSL2}), (\ref{eq:J*J}) and Lemma \ref{lem:equivofmonoidalcat}.
\epf

The following result covers the first two statements from Theorem \ref{thm:jsheaves}:
\begin{theorem}\label{thm:jsheavesa}
Let $w\in W$ and let $w=s_1\cdots s_l$ be a reduced expression of $w$ as a product of simple reflections. Then \\
(i) The object $\Kcal_w:=\Kcal_{s_1}\ast\cdots\ast\Kcal_{s_l}\in \eUop\D(G)\eUop$ is independent of the choice of the reduced expression. Note that $\Kcal_{1_W}:=\eUop$.\\
(ii) The object $\Kcal_w[-\dim \Uop]$ is an irreducible perverse sheaf which is the middle extension of a (shifted) local system supported on a closed subset of $\Bop w\Bop$.\\
(iii) If $w_1,w_2\in W$ are such that $l(w_1w_2)=l(w_1)+l(w_2)$, then we have a canonical isomorphism $\Kcal_{w_1w_2}\cong \Kcal_{w_1}\ast\Kcal_{w_2}$.
\end{theorem}

\subsection{Proof of Lemma \ref{lem:bimodule}} Let us begin by proving Lemma \ref{lem:bimodule}. Since $\eUop,\eL\in \D(G)$ are closed idempotents, it follows that the full subcategory $\eUop\D(G)\eL$ is a left $\eUop\D(G)\eUop$-module category and a right $\Wh$-module category and that these two actions are compatible. 

Let us now see on which double cosets in $\Uop\backslash G/U$ can the objects of $\eUop\D(G)\eL$ be supported. Note that from the Bruhat decomposition we obtain $G=\Uop NU$, where we recall that $N$ is the normalizer of the maximal torus $T$. By a standard argument from Mackey theory, objects of $\eUop\D(G)\eL$ can only be supported on those $g\in G$ such that $\eUop\ast\delta_g\ast\eL\neq 0$ (see \cite{BoDr:14, Desh:16}), i.e. those $g$ such that $\eUop\ast {}^g\eL\neq 0$. By \cite[Lemma 2.14]{Desh:16} this is equivalent to the condition that $\Qlcl|_{\Uop\cap { }^gU}\cong { }^g\L|_{\Uop\cap { }^gU}$. (Note that the intersection $\Uop\cap { }^gU$ is connected.) Without loss of generality let us suppose that $g\in N$ since $N$ represents all the $\Uop-U$-double cosets. But since $\L$ is a non-degenerate multiplicative local system on $U$, the isomorphism $\Qlcl|_{\Uop\cap { }^gU}\cong { }^g\L|_{\Uop\cap { }^gU}$ with $g\in N$ is possible only if $\Uop\cap{}^gU$ contains no simple root subgroup. Hence we must necessarily have $g\in T$. This proves that the objects of $\eUop\D(G)\eL$ can only be supported on $\Uop TU$ which is an open subset of $G$.

Note that the open subvariety $\Uop TU\subset G$ is isomorphic to the product $\Uop\times T\times U$. Hence it is clear that we have a triangulated equivalence $\D(T)\rar\cong\eUop\D(\Uop TU)\eL=\eUop\D(G)\eL$ given by 
\beq
\D(T)\ni C\mapsto \eUop\ast C\ast \eL.
\eeq
The inverse functor is given by
\beq
\eUop\D(G)\eL\ni C\mapsto C|_T[-4\dim U](-2\dim U).
\eeq
The other equivalences $$\D(T)\cong \eUop\D(\Bop)\cong \D(B)\eL\cong  \eUop\D(G)\eL$$ stated in Equation (\ref{eq:triangulatedequivalenceofbimodule}) are also clear.

The full subcategory $\eUop\D(\Bop)\eUop=\eUop\D(\Bop)\subset \eUop\D(G)\eUop$ is equivalent to $\D(T)$ as a triangulated monoidal category. It is clear that the left action of this subcategory on the bimodule category $\D(T)\cong\eUop\D(\Bop)\cong\eUop\D(G)\eL$ is the usual left convolution with compact support. 

From the above discussion we obtain:
\blem\label{lem:leftaction}
The left action of an object $A\in \Yok$ on an object $C\in\D(T)$ is given by 
\[
A\star (\cdot):C\mapsto (A\ast\eUop\ast C\ast \eL)|_T[-4\dim U](-2\dim U)=(A\ast C\ast \eL)|_T[-4\dim U](-2\dim U),
\]
where we have denoted the above action by $\star$, i.e. 
\beq A\star C:=(A\ast C\ast \eL)|_T[-4\dim U](-2\dim U).\eeq
\elem

\subsection{A quotient of the Yokonuma-Hecke category}
In this subsection we describe a certain quotient $\Yscr$ of the Yokonuma-Hecke category modulo a thick triangulated two-sided ideal $\Iscr\subset \eUop\D(G)\eUop$. We will see that the ideal $\Iscr$ is such that its left action on the bimodule category $\D(T)$ is trivial and hence $\D(T)$ gets the structure of a $\Yscr-\Wh$-bimodule category.

By Lemmas \ref{lem:equivofmonoidalcat}, \ref{lem:qiesalpha}, for each $\alpha\in \Delta$ we have the quasi-idempotent $e_{s_{\alpha}}=e_{\Uop_{s_\alpha}}\ast {f_{\alpha}}_!e_{SL_2}\in \Yok$, where $f_\alpha:SL_2\onto G'_\alpha\subset G_\alpha\subset G$ is the root homomorphism and $e_{SL_2}$ is the quasi-idempotent in $e_{U^-}\D(SL_2)e_{U^-}$ defined in (\ref{eq:esl2}).

\blem\label{lem:trivialleftaction}
For each $\alpha\in \Delta$ the quasi-idempotent $e_{s_\alpha}\in \Yok$ acts trivially on the left on the bimodule $\eUop\D(G)\eL$:
\[e_{s_\alpha}\ast A=0 \mbox{ for every object }A\in \eUop\D(G)\eL.\]
\elem
\bpf
Consider the ``constant'' idempotent $e_{U^+}\in \D(SL_2)$ supported on the  subgroup $U^+\subset SL_2$ of upper unipotent matrices. Since $e_{SL_2}$ is constant on all of $SL_2$, and in particular $U^+$-biequivariant for right and left translations by $U^+\cong\Ga$, we have $e_{SL_2}\cong e_{SL_2}\ast e_{U^+}\cong e_{U^+}\ast e_{SL_2}$. Since the functor ${f_\alpha}_!:\D(SL_2)\to \D(G)$ is monoidal, we obtain that $${f_\alpha}_! e_{SL_2}\cong ({f_\alpha}_! e_{SL_2})\ast ({f_\alpha}_! e_{U^+})\cong {f_\alpha}_! e_{SL_2}\ast e_{U_\alpha}\cong e_{U_\alpha}\ast {f_\alpha}_! e_{SL_2},$$ where $e_{U_\alpha}= {f_\alpha}_! e_{U^+}$ is the constant idempotent supported on the root subgroup $U_\alpha$. Hence we obtain that 
\beq\label{eq:eualpha} e_{U_\alpha}\ast e_{s_\alpha}\cong e_{s_\alpha}\ast e_{U_\alpha}\cong e_{s_\alpha}.\eeq
By Lemma \ref{lem:bimodule}, the convolution $e_{s_{\alpha}}\ast A\in \eUop\D(G)\eL=\eUop\D(\Uop TU)\eL$ can only be supported on those $\Uop-U$ double cosets represented by $T$. On the other hand by (\ref{eq:eualpha}), $e_{s_{\alpha}}\ast A\in e_{U_\alpha}\D(G)\eL$. Now for each $t\in T$,  $e_{U_\alpha}\ast\delta_t\ast\eL=0$. Hence $e_{s_{\alpha}}\ast A$ must vanish on $T$. Hence we conclude that $e_{s_{\alpha}}\ast A$ must be zero.
%By Lemma \ref{lem:equivofmonoidalcat} we have  $e_{U_{-\alpha}}\D(G_\alpha)e_{U_{-\alpha}}\cong \eUop\D(\Pop_\alpha)\eUop\subset \Yok$ with $e_{s_\alpha}\in\eUop\D(\Pop_\alpha)\eUop$ corresponding to ${f_\alpha}_!e_{SL_2}\in e_{U_{-\alpha}}\D(G_\alpha)e_{U_{-\alpha}}$. We also have the equivalences $e_{U_{-\alpha}}\D(G_\alpha)e_{\L|_{U_\alpha}}\cong\D(T)\cong \eUop\D(G)\eL$ and the left actions of                  $e_{U_{-\alpha}}\D(G_\alpha)e_{U_{-\alpha}}\cong \eUop\D(\Pop_\alpha)\eUop$ are compatible along these equivalences. Hence in order to prove the lemma we may suppose that $G=G_\alpha$.
\epf

With this lemma proved, we make the following definition:
\bdefn
Let $\Iscr\subset\Yok$ be the thick triangulated two-sided ideal generated by all the quasi-idempotent $\{e_{s_\alpha}|\alpha\in \Delta\}$. Define $\Yscr:=\Yok/\Iscr$. The  category $\D(T)\cong\eUop\D(G)\eL$ has the structure of a $\Yscr-\Wh$-module category.
\edefn

Let $\underline{W}$ denote the pointed monoidal category whose objects are the elements of $W$, the only morphisms are isomorphisms and tensor product is given by multiplication in $W$.
\bprop
The functor $\underline{W}\to \Yscr$ defined by $w\mapsto \Kcal_w$ (viewed as an object in the quotient modulo $\Iscr$) has a natural monoidal structure.
\eprop
\bpf
This follows from Theorem \ref{thm:jsheavesa} along with Lemma \ref{lem:qiesalpha}(ii) and the definition of the ideal $\Iscr$.
\epf

\subsection{Proof of Theorem \ref{thm:jsheaves}(iii)}
We now compute the left action of the objects $\Kcal_w$ on the bimodule category $\D(T)$ and complete the proof of Theorem  \ref{thm:jsheaves}(iii). It is enough to consider the left action of $\Kcal_{s_\alpha}$ on $\D(T)$ for a simple root $\alpha\in \Delta$. Let $C\in \D(T)$. By Lemmas \ref{lem:equivofmonoidalcat} and \ref{lem:leftaction} the left action is given by 
\beq
\Kcal_{s_\alpha}\star C = (e_{\Uop_{s_\alpha}}\ast {f_{\alpha}}_!\Kcal \ast C\ast \eL)|_T[-4\dim U](-2\dim U).
\eeq
Note that we have the identifications $\eUop\D(G)\eL\cong \D(T)\cong e_{U_{-\alpha}}\D(G_\alpha)e_{\L_\psi}$, where $\L_\psi$ is the Artin-Schreier multiplicative local system on $U_\alpha\xto{\cong}\Ga$. Using this we can see that the action above can be computed within the subgroup $G_\alpha\subset G$ which is of semisimple rank 1:
\beq
\Kcal_{s_\alpha}\star C= ({f_{\alpha}}_!\Kcal \ast C\ast e_{\L_\psi})|_T[-4](-2).
\eeq
Let us fix $g\in T$ and compute the stalk of $\Kcal_{s_\alpha}\star C$ at $g$:
$${f_{\alpha}}_!\Kcal \ast C\ast e_{\L_\psi}[-4](-2)(g)=\int\limits_{\substack{{x\in SL_2,t\in T,a\in U_\alpha}\\{f_\alpha(x)ta=g}}}\Kcal(x)\otimes C(t)\otimes \L_\psi(a)[-2](-1)$$
$$=\int\limits_{\substack{{x=ru\in B^+,t\in T,a\in U_\alpha}\\{f_\alpha(x)ta=g}}}\Kcal(x)\otimes C(t)\otimes \L_\psi(a)[-2](-1)\ \  \mbox{ (where $B^+=\left\{\begin{tiny}\begin{pmatrix}
    r\hspace{-7pt} & ru\\
    0\hspace{-7pt} & r^{-1}
\end{pmatrix}=\begin{pmatrix}
    r\hspace{-7pt} & 0\\
    0\hspace{-7pt} & r^{-1}
\end{pmatrix}\begin{pmatrix}
    1\hspace{-7pt} & u\\
    0\hspace{-7pt} & 1
\end{pmatrix}\end{tiny}
\right\}\subset SL_2$)}$$
$$=\int\limits_{\substack{{r\in \Gm, u\in \Ga,t\in T,a\in U_\alpha}\\{\alpha^\vee(r)f_\alpha(u)ta=g}}}\L_\psi(ru)\otimes C(t)\otimes \L_\psi(a)[+2](+1) \ \ \mbox{ (using the definition of $\Kcal$)}$$
$$=\int\limits_{\substack{{r\in \Gm, a,u\in \Ga,t\in T}\\{\alpha^\vee(r)t=g}\\\alpha(t^{-1})u+a=0}}\L_\psi(ru+a)\otimes C(t)[2](1)$$
$$=\int\limits_{r\in \Gm, u\in \Ga}\L_\psi\left(\left(1-\alpha(g^{-1})r\right)ru\right)\otimes C(\alpha^\vee(r^{-1})g)[2](1)\ \ \mbox{ (substituting $t,a$ and using $\alpha(\alpha^\vee(r))=r^2$)}$$
$$=\int\limits_{\substack{r=\alpha(g)\\u\in \Ga}}C(\alpha^\vee(r^{-1})g)[2](1)\ \ \mbox{ (using cohomology vanishing of non-trivial multiplicative local systems on }\Ga)$$
$$=\int\limits_{u\in \Ga}C(\alpha^\vee(\alpha(g^{-1}))g)[2](1)=C(s_\alpha(g)).$$
Hence we conclude that the left action of $\Kcal_{s_\alpha}$ on $\D(T)$ coincides with the natural action of $s_\alpha\in W$ on $\D(T)$. This completes the proof of Theorem \ref{thm:jsheaves}(iii).

\section{Bi-whittaker sheaves and central sheaves}\label{sec:ptm2}
In this section we will complete the proof of Theorem \ref{thm:main2}
\subsection{The functor $\Wh\to \Dcirc_W(T)$}\label{sec:whtocentral}
In this subsection, we construct the triangulated monoidal functor from the bi-Whittaker category to the category of central sheaves on the torus using the structure of the $\Yscr-\Wh$-bimodule category on $\D(T)$. 

We have seen that we have a monoidal functor $\D(T)\cong\eUop\D(\Bop)\subset \Yok\to \Yscr$ and that the induced left action of $\D(T)$ on the $\Yscr-\Wh$-bimodule $\D(T)$ is the usual convolution. Hence any left $\Yscr$-module category endofunctor of $\D(T)$ is a pair $(C,\phi)$ where $C\in \D(T)$ and $\phi$ is a family of functorial isomorphisms for $A\in \Yscr, D\in \D(T)$:
$$\phi_{A,D}:(A\star D)\ast C\xto{\cong}A\star(D\ast C)$$
satisfying certain compatibilities. For example using the compatibilities and the monoidal functor $\D(T)\to \Yscr$ the natural isomorphisms above are determined by the isomorphisms \beq\phi_A:(A\star \delta_1)\ast C\xto{\cong}A\star C.\eeq
In particular since we have a monoidal functor $\underline{W}\to \Yscr$, for each $w\in W$ we have the isomorphisms (taking $A=\Kcal_w$)
\beq\label{eq:wequivonbiwhit}\phi_w:C\xto{\cong}\Kcal_w\star C=w(C)\eeq
and hence the pair $(C,\phi)$ induces a $W$-equivariant structure on $C$.

Now the right action of an object of $\Wh$ on the $\Yscr-\Wh$-bimodule $\D(T)$ induces a $\Yscr$-module endofunctor of $\D(T)$. By the above observations, this defines a monoidal functor $\Wh\to \D_W(T)^\op$. But $\D_W(T)$ is a symmetric monoidal category, hence we obtain a monoidal functor
\beq
\Wh\to \D_W(T).
\eeq
Note that this is a triangulated functor since the composition of the above with the forgetful functor to $\D(T)$ is given by $\Wh\ni A\mapsto (\eUop\ast A)|_T[-4\dim U](-2\dim U)\in \D(T)$, which is a triangulated functor.
\bprop
There is a natural monoidal functor \[\Ecal nd_\Yscr(\D(T))\to \D_W(T)\]
whose image lies inside the subcategory $\Dcirc_W(T)$ of central sheaves on the torus. In particular we have a natural monoidal functor 
\beq
\xi:\Wh\to \Dcirc_W(T).
\eeq
Under this functor, an object $A\in \Wh$ maps to the pair $\left((\eUop\ast A)|_{T}[-4\dim U](-2\dim U),\rho^A\right)$, where $\rho^A$ is the $W$-equivariant structure coming from the left $\Yscr$-module endofunctor structure of the right action of $A$ on the bimodule category $\D(T)$.
\eprop
\bpf
Let $(C,\phi)\in \Ecal nd_\Yscr(\D(T))$ in the previous notation. Let $\alpha\in \Delta$ and $\Kcal_{s_\alpha}\in \Yok$ the corresponding object which is supported on $\Bop\sqcup\Bop s_\alpha \Bop$. Let us look at the restriction of this object to the subcategory $\eUop\D(\Bop)\cong \D(T)$. In the case of $G=SL_2$, this restriction corresponds to the (shifted) constant sheaf on the maximal torus $\Gm\subset SL_2$. Hence we see that $\Kcal_{s_\alpha}|_{\Bop}$ corresponds to the object $\alpha^\vee_!\Qlcl[2](1)\in \D(T)$. We have the canonical adjunction morphism in $\Yok$
\[\Kcal_{s_\alpha}\to {\Kcal_{s_\alpha}}_{\Bop}=\eUop\ast\alpha^\vee_!\Qlcl[2](1).\]
Let $\chi\in \Ccal(T)$ be a Kummer local system such that ${\alpha^\vee}^*\chi$ is the trivial Kummer local system on $\Gm$. Let $\chi^\vee$ be the dual Kummer local system. Then we also have a canonical isomorphism 
\[{\alpha^\vee}^*\chi^\vee\xto{\cong}\Qlcl\]
of Kummer local systems.
By adjunction and then applying Verdier duality, we get the morphisms
\[\chi^\vee\to \alpha^\vee_*\Qlcl\mbox{ and }\alpha^\vee_!\Qlcl[2](1)\to \chi[2\dim T](\dim T)=:e_\chi.\]
Here $e_\chi=\chi[2\dim T](\dim T)$ denotes the minimal quasi-idempotent in $\D(T)$ corresponding to the Kummer local system $\chi\in \Ccal(T)$.
This gives us the composition
\beq
\Kcal_{s_\alpha}\to {\Kcal_{s_\alpha}}_{\Bop}=\eUop\ast\alpha^\vee_!\Qlcl[2](1)\to \eUop\ast e_\chi.
\eeq
Using the object $(C,\phi)$ we get the commutative diagram
\beq
\xymatrix{
C\ar[r]^-{\phi_{s_\alpha}}\ar[d] & \Kcal_{s_\alpha}\star C=s_\alpha(C)\ar[d]\\ 
e_\chi\ast C\ar@{=}[r] & e_\chi\ast C.}
\eeq
An easy computation shows that $e_\chi\ast C=H^*_c(T,C\otimes \chi^\vee)\otimes e_\chi$. Hence the above commutative diagram means that the action of $s_\alpha$ on $H^*_c(T,C\otimes\chi^\vee)$ induced by the $W$-equivariant structure on $C$ is trivial. Hence by Lemma \ref{lem:centralconditions}, the $W$-equivariance structure on $C$ coming from the pair $(C,\phi)$ is central, completing the proof of the proposition.
\epf

\subsection{Equivalence of two $W$-equivariance structures}\label{sec:equivalenceofWequiv}
Composing $\xi$ with the central monoidal functor $\HC_\L$ we obtain $\D_G(G)\xto{\HC_\L}\Wh\xto{\xi} \Dcirc_W(T)$. The right action of an object $A\in \Dcirc_G(G)$ on $\D(T)\cong \eUop\D(G)\eL$ maps the object $\eUop\ast\eL$ to $\eUop\ast A\ast \eL$. Since $A\in \Dcirc_G(G)$, $\eUop\ast A$ is supported on $\Bop$ and hence the object $\eUop\ast A\ast \eL$ corresponds to the object $\res(A)\in \D(T)$. For $w\in W$, the $w$-equivariance structure on $\xi(\HC_\L(A))$ corresponds to an isomorphism $\eUop\ast A\ast\eL\cong \Kcal_w\ast\eUop\ast A\ast\eL$. This isomorphism comes from the the natural braiding isomorphism $\eUop\ast A\ast \Kcal_w\xto{\cong} \Kcal_w\ast\eUop\ast A$ convolved with the idempotent $\eL$ on the right.

On the other hand for an object $A\in \Dcirc_G(G)$, we have the object $\res^\circ(A)\in \Dcirc_N(T)$ and forgetting the $T$-equivariance, an object of $\Dcirc_W(T)$. To complete the proof of Theorem \ref{thm:main2} we will need the following:
\bprop\label{prop:naturalequivofcompositions}
The following two compositions of functors are naturally equivalent
\beq\label{eq:DGtoDW1}\Dcirc_G(G)\xto{\res^\circ}\Dcirc_N(T)\rar{}\Dcirc_W(T),\eeq
\beq\label{eq:DGtoDW2}\Dcirc_G(G)\xto{\HC_\L}\Wh\rar{\xi}\Dcirc_W(T).\eeq
\eprop
\bpf
Let $A\in\Dcirc_G(G)$. The composition of the two functors above with the forgetful functor $\Dcirc_W(T)\to \D(T)$ map $A$ to $\eUop\ast A\in \eUop\D(\Bop)\cong \D(T)$ which corresponds to the object $\res(A)\in \D(T)$. Hence we only need to verify that the two above functors define the same $W$-equivariance structure on $\res(A)$. For this it suffices to compare the two $s_\alpha$-equivariance structures on $\res(A)$ for each simple root $\alpha$. By parabolic restriction and using Lemma \ref{lem:equivofmonoidalcat} we are reduced to the case where $G$ is of semisimple rank 1. Let $W=\{1,s\}$ be the Weyl group of $G$ in this case.

Let us first consider the case where $G\cong SL_2\times S$, where $S$ is a torus. Then $\Uop\backslash(SL_2\times S)\cong (\mathbb{A}^2\setminus 0)\times S\subset \bbA^2\times S$. Let $M_2$ be the algebraic (multiplicative) monoid of $2\times 2$-matrices and let $\bop\subset M_2$ be the submonoid formed by all lower triangular matrices $\begin{pmatrix}t & 0\\  c & d\end{pmatrix}$ with $t\in \Gm$. Note that $\Uop\backslash\bop\cong \Gm\times\bbA^1$ as algebraic monoids. Consider the right multiplication action of the monoid $\bop\times S$ on $\bbA^2\times S$. Also note that $\D(\bop\times S)$ is a monoidal category under convolution with compact supports and that $\eUop\D(\bop\times S)=\eUop\D(\bop\times S)\eUop$ is a monoidal category (with unit $\eUop$) which is monoidally equivalent to $\D(\Gm\times \bbA^1\times S)$.

On the other hand, the left action of $T\cong\Uop\backslash\Bop\cong\Gm\times S$ on $\Uop\backslash G\cong (\bbA^2\setminus 0)\times S$ is the natural one and extends to an action on $\bbA^2\times S$. We see that the left action of $\Gm\times S$ on $\bbA^2\times S$ commutes with the right action of $M_2\times S$ and hence with that of $\bop\times S$. For both the left and right actions $0\times S\subset \bbA^2\times S$ is a closed invariant subset.

From the above it follows that $\eUop\D(G)\cong \D(\Uop\backslash G)\cong \D((\bbA^2\setminus 0)\times S)$ is a $\eUop\D(\Bop)-\D(\bop\times S)$-bimodule (considered as a quotient of the bimodule $\D(\bbA^2\times S)$). Hence we have that the Hecke category $\eUop\D(G)\eUop$ is a $\eUop\D(\Bop)-\eUop\D(\bop\times S)$-bimodule category. Consider the closed embedding $\iota:\bbA^1\xhookrightarrow{}\Gm\times \bbA^1\times S$ defined by $x\mapsto (1,x,1)$ and the object $\iota_!\L\in \D(\Gm\times\bbA^1\times S)$. Then we have 
\beq\label{eq:deltaandKL}\eUop\ast\delta_{s}\ast \eUop\ast \iota_!\L\cong  \Kcal_s, 
\eeq
where $s=\left(\begin{pmatrix}
    0 & 1\\
    -1 & 0
\end{pmatrix},1\right)\in N_G(T).$

We have the natural inclusions $\Bop=B^-\times S\subset \bop\times S,\  T=\Gm\times S\xhookrightarrow{(t,x)\mapsto(t,t^{-1},x)} \Gm\times \bbA^1\times S$ and hence a symmetric monoidal functor $\eUop\D(\Bop)\cong \D(T)\to \D(\Gm\times \bbA^1\times S)\cong \eUop\D(\bop\times S)$. We also obtain the sequence of braided monoidal functors $\Dcirc_G(G)\to \eUop\D_{\Bop}(\Bop)\to \eUop\D(\Bop)\cong \D(T)\to \D(\Gm\times \bbA^1\times S).$

For $A\in \Dcirc_G(G)$ the two $s$-equivariance structures on $\res(A)\in \D(T)$ in the proposition come from the braiding isomorphism of the object $\eUop\ast A$ with the objects $\eUop\ast\delta_s\ast\eUop$ and $\Kcal_s$ in $\eUop\D(G)\eUop$ respectively, using the central monoidal structure of the functor $\Dcirc_G(G)\to \eUop\D(G)\eUop$. The equality of the two $s$-equivariance structures then follows from (\ref{eq:deltaandKL}) completing the proof for $G=SL_2\times S$.

A general reductive group $G'$ of semisimple rank 1 can be realized as a quotient $G'=G/\Lambda$, where $G$ is as before and $\Lambda$ is a finite central subgroup of $G$. Then for $A'\in \Dcirc_{G'}(G')$ it follows that the two $s$-equivariance structures on $\res(A')\in \D(T')$ are equal using the result for $G$, proper base change and the fact that the pullback $\D(T')\to \D(T)$ is a faithful functor.
\epf

\subsection{Completion of the proof}
In Section \ref{sec:whtocentral} we have constructed the monoidal functor 
\[\xi:\Wh\to \Dcirc_W(T).\]
Going in the other direction, we have the monoidal functor (see Corollary \ref{cor:wactiononindcirc})
\[\Dcirc_W(T)\xto{}\Dcirc_N(T)\xto[\cong]{{\ind^\circ}^W}\Dcirc_G(G)\xto{\HC_\L}\Wh.\]

By Lemma \ref{prop:naturalequivofcompositions}, the composition
\beq
\Dcirc_W(T)\to\Dcirc_N(T)\xto{\cong}\Dcirc_G(G)\xto{\HC_\L}\Wh\rar{\xi}\Dcirc_W(T)\eeq
is naturally equivalent to the composition
\beq\label{eq:DcirctoDcirc}
\Dcirc_W(T)\to\Dcirc_N(T)\xto{\cong}\Dcirc_G(G)\xto{\res^\circ}\Dcirc_N(T)\rar{}\Dcirc_W(T).\eeq
The above composition is naturally equivalent to the identity functor on $\Dcirc_W(T)$.

Next we consider the opposite composition
\beq\Wh\xto{\xi}\Dcirc_W(T)\to \Dcirc_N(T)\xto{\cong}\Dcirc_G(G)\xto{\HC_\L}\Wh.\eeq

\blem\label{lem:WhtoWh}
(i) The composition
\beq\label{eq:HecketoWh}
\eUop\D_{\Bop}(G)\xto{\av_{\Bop}^G}\D_G(G)\xto{\HC_\L}\Wh
\eeq
is given by 
\beq
\eUop\D_{\Bop}(G)\ni A\mapsto \L\ast A\ast\eL.
\eeq
(ii) The composition
\beq\label{eq:WhtoWh}\Wh\xto{\xi}\Dcirc_W(T)\to\Dcirc_N(T)\to \D_T(T)\hookrightarrow\eUop\D_{\Bop}(G)\xto{\av_{\Bop}^G}\D_G(G)\xto{\HC_\L}\Wh\eeq
is given by convolution with the object $\L\ast\eUop\ast\eL\in \Wh$.
\elem
\bpf
Let us first consider the composition (\ref{eq:HecketoWh}). Under this composition an object $A\in \eUop\D_{\Bop}(G)$ maps to $\av_{\Bop}^GA \ast \eL$. Computing the stalk at $x\in G$ we have
\[\av_{\Bop}^GA \ast \eL(x)=\int\limits_{u\in U}\av_{\Bop}^GA(xu^{-1})\otimes \eL(u)\]
\[=\int\limits_{{u\in U}}\ \ \int\limits_{g\Bop\in G/\Bop}{ }^gA(xu^{-1})\otimes \eL(u)\]
\[=\int\limits_{g\Bop\in G/\Bop}\ \ \int\limits_{u\in U}{ }^gA(xu^{-1})\otimes e_{\L}(u)\]
\beq\label{eq:tobecont}=\int\limits_{g\Bop\in G/\Bop}{}^gA\ast e_{\L}(x).\eeq
Note that we have the Bruhat decomposition $G=UW\Bop$. Suppose that the coset $g\Bop\in G/\Bop-U\Bop/\Bop$. Then we claim that ${}^gA\ast\eL=0$. Suppose that $g=un$ with $n\in N-T$. Then ${}^n\Uop$ and $U$ both contain a simple root subgroup $U_\alpha$ for some $\alpha\in\Delta$. By definition of non-degeneracy, $\L$ is non-trivial on each simple root subgroup. Hence we conclude that $e_{{}^n\Uop}\ast\eL=0$, and hence $e_{{}^{un}\Uop}\ast e_{{}^u\L}=e_{{}^g\Uop}\ast\eL=0$. Now since $A\in \eUop\D_{\Bop}(G)$, we have $A\ast\eUop\cong A$. Hence ${}^gA\ast e_{{}^g\Uop}\cong{}^gA$. From the above statements we conclude that ${}^gA\ast\eL=0$ for $g\in G-U\Bop$.

Now continuing from (\ref{eq:tobecont}) we get
\[\av_{\Bop}^GA \ast \eL(x)=\int\limits_{u\in U}{}^uA\ast e_{\L}(x)
\]
\[=\int\limits_{u\in U}{}^u(A\ast e_{\L})(x)=\int\limits_{u\in U}A\ast e_{\L}(u^{-1}xu)\]
\[=\int\limits_{u\in U}A\ast e_{\L}(u^{-1}x)\otimes\L(u)=\int\limits_{u\in U}\L(u)\otimes (A\ast e_{\L})(u^{-1}x)\]
\[=\L\ast A\ast \eL(x).\] This completes the proof of statement (i).

To prove (ii) consider an object $A\in \Wh$. Under the composition (\ref{eq:WhtoWh}), $A$ gets mapped to the following object in $\Wh$:
\[\av_{\Bop}^G((\eUop\ast A)|_{\Bop})\ast\eL[-2\dim U](-\dim U)=\av_{\Bop}^G(\eUop\ast A|_{\Bop})\ast\L.\]
By part (i), this equals 
\[\L\ast\eUop\ast A_{\Bop}\ast\eL[-2\dim U](-\dim U)=\L\ast(\eUop\ast A)_{\Bop}\ast\eL[-2\dim U](-\dim U).\]
Now since $A\in\Wh$, $\eUop\ast A\in \eUop\D(G)\ast\eL$ and hence is supported on $\Bop U$. Hence we have 
\[(\eUop\ast A)_{\Bop}\ast\eL[-2\dim U](-\dim U)=\eUop\ast A.\]
Hence we see that under the composition (\ref{eq:WhtoWh}), $A\mapsto \L\ast\eUop\ast A=\L\ast\eUop\ast\eL\ast A$, proving statement (ii).
\epf

Finally we complete the proof of Theorem  \ref{thm:main2} below:
\bcor\label{cor:Waction}
(i) We have $\eL\ast\Spr=\L\ast\eUop\ast\eL$ is a semisimple perverse sheaf (shifted by $\dim U$) equipped with a $W$-action, where $\Spr=\av_{\Bop}^G\eUop$ is the Springer sheaf.\\
(ii) We have $(\L\ast\eUop\ast\eL)|_B=\eL$ and that this $\eL$ is the direct summand of $\L\ast\eUop\ast\eL$ corresponding to taking the $W$-invariants.\\
(iii) The composition
\beq\label{eq:WinvWhtoWh}\Wh\xto{\xi}\Dcirc_W(T)\to\Dcirc_N(T)\xto{{\ind^\circ}^W}\Dcirc_G(G)\xto{\HC_\L}\Wh\eeq
is naturally equivalent to the identity functor. In particular, $\xi:\Wh\to \Dcirc_W(T)$ is a triangulated monoidal equivalence with inverse $\Dcirc_W(T)\to \Dcirc_N(T)\xto{{\ind^\circ}^W}\Dcirc_G(G)\xto{\HC_{\L}}\Wh$.
\ecor
\bpf
That $\eL\ast\Spr=\L\ast\eUop\ast\eL$ follows from Lemma \ref{lem:WhtoWh}(i) and the definition of the Springer sheaf. That it is semisimple perverse (shifted by $\dim U$) follows from the results of \cite{BBM:04,BBM22}. Finally, the $W$-action on the Springer sheaf $\Spr$ gives us the desired $W$-action on $\L\ast\eUop\ast\eL$, proving (i).

That $(\L\ast\eUop\ast\eL)|_B=\eL$ is a straightforward computation. Recall that we have normalized the $W$-action on $\Spr$ so that we have $\Spr^W=\delta_1\in\D_G(G)$. Hence the $W$-invariant direct summand of $\L\ast\eUop\ast\eL$ equals $\HC_\L(\delta_1)=\eL$, proving (ii).

Finally, we note that the composition (\ref{eq:WhtoWh}) is equivalent to
\beq
\Wh\xto{\xi}\Dcirc_W(T)\to\Dcirc_N(T)\xto{\ind^\circ} \D_G(G)\xto{\HC_\L}\Wh
\eeq
and that the functor $\ind^\circ$ has a $W$-action coming from the $W$-action on $\Spr$ (see Corollary \ref{cor:wactiononindcirc}). Taking the $W$-invariants, we get the composition (\ref{eq:WinvWhtoWh}). By (ii) this corresponds to convolution with the $W$-invariant summand $\eL$ of $\L\ast\eUop\ast\eL$. Since $\eL$ is the unit object in $\Wh$, the composition (\ref{eq:WinvWhtoWh}) is naturally equivalent to the identity. On the other hand we have already seen that the composition (\ref{eq:DcirctoDcirc}) is naturally equivalent to the identity. This completes the proof of statement (iii) and hence also of Theorem \ref{thm:main2}.
\epf

\bibliographystyle{alpha}
\bibliography{papers}
\end{document}